\documentclass[review,12pt,authoryear]{elsarticle}

\usepackage{booktabs}
\usepackage{float}
\usepackage{lineno}

\usepackage{rotating}
\makeatletter

\makeatletter
\newcommand{\Rmnum}[1]{\expandafter\@slowromancap\romannumeral #1@}
\makeatother

\makeatletter 

\usepackage{ulem}
\usepackage{eurosym}
\usepackage{caption}
\usepackage{subcaption}

\usepackage{hyperref}

\usepackage{amsmath}
\usepackage{bm}  
\usepackage{cleveref}
\usepackage{ dsfont }
\usepackage{array}
\usepackage{amssymb}
\usepackage{float}
\usepackage{xspace}
\usepackage{amsthm}
\usepackage{algorithmic}
\usepackage{algorithm}
\usepackage{graphicx}
\usepackage{pifont}

\usepackage[latin1]{inputenc}
\usepackage{tikz}
\usepackage{setspace}
\usetikzlibrary{shapes,arrows}

\usepackage[authoryear]{natbib}
\usepackage{graphicx}
\usepackage{stfloats}
\usepackage{amssymb}
\usepackage{amsmath}
\usepackage{float}
\usepackage{amssymb}
\usetikzlibrary{mindmap}
\usepackage{verbatim}
\usepackage{verbatim}
\usepackage{booktabs}
\usepackage{multirow}

\usepackage{rotating}
\usepackage{makecell}

\usepackage{colortbl}
\usepackage{xcolor} 
\definecolor{blue}{RGB}{0, 0, 0} 
\usepackage{upgreek}

\makeatletter
\def\ps@pprintTitle{%
   \let\@oddhead\@empty
   \let\@evenhead\@empty
   \let\@oddfoot\@empty
   \let\@evenfoot\@empty
}
\makeatother

\begin{document}
\begin{frontmatter}
\title{Designing a Robust and Cost-Efficient Electrified Bus Network with Sparse Energy Consumption Data}


\author[1]{Sara Momen\texorpdfstring{\corref{cor1}}{}}
\author[2]{Yousef Maknoon}
\author[1]{Bart van Arem}
\author[1]{Shadi Sharif Azadeh}
\cortext[cor1]{Corresponding author. Email: s.sharifazadeh@tudelft.nl}
\address[1]{Transport and Planning Department, Delft University of Technology, Netherlands}
\address[2]{Engineering, Systems and Services Department, Delft University of Technology, Netherlands}
\begin{abstract}
This paper addresses the challenges of charging infrastructure design (CID) for electrified public transport networks using Battery Electric Buses (BEBs) under conditions of sparse energy consumption data. Accurate energy consumption estimation is critical for cost-effective and reliable electrification but often requires costly field experiments, resulting in limited data. To address this issue, we propose two mathematical models designed to handle uncertainty and data sparsity in energy consumption. The first is a robust optimization model with box uncertainty, addressing variability in energy consumption. The second is a data-driven distributionally robust optimization model that leverages observed data to provide more flexible and informed solutions.
To evaluate these models, we apply them to the Rotterdam bus network. Our analysis reveals three key insights: (1) Ignoring variations in energy consumption can result in operational unreliability, with up to 55\% of scenarios leading to infeasible trips. (2) Designing infrastructure based on worst-case energy consumption increases costs by 67\% compared to using average estimates. (3) The data-driven distributionally robust optimization model reduces costs by 28\% compared to the box uncertainty model while maintaining reliability, especially in scenarios where extreme energy consumption values are rare and data exhibit skewness. In addition to cost savings, this approach provides robust protection against uncertainty, ensuring reliable operation under diverse conditions.
\end{abstract}
\end{frontmatter}

\section{Introduction}

\textcolor{blue}{Countries worldwide are increasingly striving to reduce their reliance on fossil fuel services. For instance, the Netherlands has made notable advancements, with The Hague taking a pioneering step by banning advertisements for fossil fuel-based services. In this paradigm, zero-emission urban mobility has emerged as a critical focus area.} Battery Electric Buses (BEBs) are leading this change due to their independence from traditional infrastructure like cables and rails. This flexibility allows BEBs to integrate seamlessly into existing transport systems. They offer a scalable solution for cities to reduce their environmental impact. However, designing effective and reliable BEB networks poses challenges. A key issue is estimating the energy consumption and strategically placing charging stations to ensure continuous service.

Accurately estimating energy consumption is critical for the cost-efficient and reliable deployment of electrified transport networks. This estimation is challenging due to a wide range of influential factors, including traffic conditions, passenger load, weather, road characteristics, and driving behavior \citep{chen2021data}. Inaccurate estimates of energy consumption impact electrification costs, battery lifetime, and service levels \citep{shadi-azadeh2022electrification}.

In the literature, several methods have been proposed to estimate energy consumption. These methods can be categorized into two main groups. The first group comprises approaches developed based on physics-based and chemical process models to determine energy consumption \citep{hjelkrem2021battery, teichert2019joint, Ricardo-scarinci2019electrification, cited-shadi_zhou2023electric}. While advanced, these approaches are limited to laboratory settings and do not account for interactions between operational, topological, and environmental factors. \textcolor{blue}{The second group comprises data-driven methods, which leverage large datasets from real-world BEB operations to estimate energy consumption \citep{gallet2018estimation, pamula2020estimation, al2020estimation}. These models excel at identifying interactions between various factors within existing systems. However, their reliance on extensive data streams restricts their applicability to well-established BEB networks. Additionally, since these models are trained on specific operational conditions, their predictions often lack generalizability to different contexts or new deployments. As a result, operators planning new BEB systems must rely on costly field tests to collect the necessary data, leading to significant uncertainty and variability in energy consumption estimates \citep{Wagenknecht2017}.}

\textcolor{blue}{This study addresses the challenges of sparse energy consumption data in designing charging infrastructure for Battery Electric Buses (BEBs), a problem known as Charging Infrastructure Design (CID). By optimizing the placement and types of charging stations (e.g., fast-charging and standard stations) and determining onboard battery capacities, the proposed approach ensures reliable and cost-effective electrification. Incorporating data sparsity and operational uncertainty, our framework overcomes limitations of existing methods and enables robust, scalable infrastructure solutions for diverse deployment scenarios.}

\textcolor{blue}{We first model the CID problem using average energy consumption estimates, acknowledging that neglecting energy consumption uncertainty can lead to designs that fail to ensure operational feasibility. To address this challenge in data-scarce environments, we propose two robust optimization models tailored to the CID problem. The first, a box uncertainty robust model with a budget (BoU-CID), uses an uncertainty box framework to account for variability by defining an average value and its maximum deviation. This approach is effective when only the range of uncertainty is known but can result in highly conservative designs. To balance cost-efficiency and reliability, the second model, a data-driven distributionally robust chance constraint approach (DRCC-CID), incorporates observed data characteristics into the robust optimization framework, producing more practical and adaptive solutions. Both models are evaluated under comprehensive scenarios to assess their performance in terms of cost-efficiency, design reliability, and computational time, demonstrating their potential for robust and scalable CID solutions.}

\textcolor{blue}{We summarize the contributions of this paper as follows: We introduce two mathematical models for designing the charging infrastructure for BEBs under the constraints of uncertain and sparse energy consumption data and investigate the impact of different approaches to modeling uncertainty on the system's total cost. The scalability and performance of the proposed models are tested on a series of test instances as well as the Rotterdam bus network. The results demonstrate that incorporating sparse energy consumption data significantly reduces total electrification costs compared to the range-based approach, even under highly conservative design scenarios. Furthermore, we evaluate the long-term costs of BEB battery lifetimes by analyzing charge-discharge cycles until battery failure across different models. Our findings show that integrating sparse data into robust design leads to a charge-discharge pattern that minimizes long-term battery-related costs. Lastly, we highlight the added value of data collection, quantifying how increasing the number of observations improves BEB charging infrastructure design in terms of computational efficiency and reliability.}

The remainder of this paper is organized as follows. Relevant literature is discussed in \Cref{sec:LR}. \textcolor{blue}{The elements of the CID problem and underlying assumptions as well as details of the mathematical models are provided in \Cref{sec:MF}. In \Cref{sec:results}, computation results and the case study are discussed.} Finally, our findings are summarized in \Cref{sec:conclusion}.

\section{Related Literature}
\label{sec:LR}

Electrifying public transport systems with BEBs involves optimizing trade-offs between the strategic placement of charging stations, the selection of charger types \textcolor{blue}{(e.g., fast-feeding and standard charging station)}, and the determination of appropriate battery sizes. The introduction of fast charging technology enables buses to recharge rapidly at stops or terminals, achieving a full recharge in less than a minute \cite{auge2015keynote}.  \textcolor{blue}{In this study, two types of charging stations are considered: fast-feeding and standard charging stations. Fast-feeding stations allow rapid recharging, capable of restoring a battery to 80\% capacity within 15 to 20 seconds \citep{auge2015keynote,tomaszewska2019lithium}. In contrast, standard stations recharge batteries at a slower rate, with the amount of charge depending on dwell time and charger power. The efficiency of electrification deployment critically hinges on accurate energy consumption estimation.} 

The optimal solution balances these factors to ensure operational reliability, which can be expressed through specific bus schedules or generalized bus frequency, while considering costs of electrification. Finding a reliable, cost-efficient solution is heavily influenced by the availability and accuracy of data on energy consumption. In this section, we provide an overview of the state-of-the-art related to the planning problem. For a comprehensive review of optimization problems for electric buses, readers are referred to \cite{zhou2024charging}.

Studies on the design of BEB networks utilize three perspectives to model energy consumption: point estimation values, simulation-based approaches to identify variations, and random variable.

The first group focuses on the design problem using nominal energy consumption values. In \cite{xylia2017locating} and \cite{shadi-azadeh2022electrification}, the authors consider nominal energy consumption values. \cite{xylia2017locating} studied the optimal placement of charging stations for a mixed fleet of biodiesel, electric, and biogas buses in Stockholm, providing an evaluation framework for deploying an alternative-fuel public transport network and performing sensitivity analysis to assess solution robustness. \cite{shadi-azadeh2022electrification} focused on battery lifetime as an additional component in the electrification process, formulating the problem as a bi-objective optimization with total ownership cost and battery lifetime as objectives. They considered nominal energy consumption values and controlled robustness by imposing fixed upper and lower Bounds on battery charge. Finally, \cite{he2022integrated} studied the BEB deployment problem, incorporating detailed charging schedules. The authors proposed a two-phase optimization framework: the first phase uses conservative (worst-case) energy consumption values for system deployment and charging schedules, while the second phase employs a rolling horizon approach to adapt the bus schedule plan.

For BEBs, energy consumption is influenced by factors such as passenger load, travel time, weather conditions, and road characteristics. The first group of studies evaluates the robustness of solutions through post-analysis of energy consumption variability using scenarios. In contrast, the second group considers energy consumption variability a priori, employing simulation models. \cite{sinhuber2012study} considered factors such as vehicle speed and route stops, proposing a simulation model to understand energy consumption variation and inform battery size decisions. \cite{rogge2015fast} further utilized this simulation model, proposing a hierarchical approach to determine fleet size, charging station placement, and battery capacity. In \cite{sebastiani2016evaluating} and \cite{teichert2019joint}, the authors integrate simulation into a simulation-optimization framework. \cite{sebastiani2016evaluating} simulate energy consumption variations based on bus speed and passenger load. This simulation is embedded in an optimization model to determine charging station placement and bus charging times. \cite{teichert2019joint} study an existing system with a fixed bus schedule and characteristics. They assume uniform environmental conditions, with energy consumption influenced only by travel time variation. Their simulation model determines energy consumption variation. A numerical approach then minimizes ownership cost by optimizing battery size and the number and type of charging stations. \textcolor{blue}{Simulation models and expert judgments are commonly used to estimate system operations, but their predictions often deviate from empirical observations, failing to capture real-world variability in energy consumption. Recognizing these limitations, a second group of studies treats energy consumption as a random variable to explicitly incorporate uncertainty into the modeling process.}

\textcolor{blue}{In the stochastic programming framework, random variables are assumed to follow known probability distributions, allowing uncertainty to be addressed using computationally intensive scenario-based methods. This approach is particularly effective when reliable probabilistic information is available. For example, \cite{Kun-an2020battery} accounted for stochastic variations in electricity demand due to weather and traffic fluctuations, using stochastic integer programming to optimize charging station placement and fleet size. Similarly, \cite{cited-shadi_zhou2023electric} modeled the state of charge as a stochastic variable influenced by travel time and battery degradation. By employing a two-stage stochastic programming approach, they optimized charging facility deployment and operational scheduling.}

\textcolor{blue}{The second approach, robust optimization, addresses uncertainty by solving problems under the worst-case realization of random variables. This method has been widely applied to BEB-related challenges, with studies incorporating diverse sources of uncertainty into their formulations. For instance, \cite{hu2022joint} accounts for travel time and passenger demand uncertainty to optimize charging station placement and bus schedules, while \cite{gairola2023optimization} focuses on energy consumption variability with similar objectives. Similarly, \cite{bai2022robust} addresses fleet sizing and wireless bus design, incorporating uncertainties in energy consumption and charging supply, and \cite{iliopoulou2021robust-BoU-partc} optimizes route planning and charging infrastructure under power supply and electricity demand uncertainty. Expanding on this, \cite{avishan2023electric} integrates energy consumption and travel time uncertainty into electric bus scheduling, minimizing costs related to scheduling, fleet purchasing, energy procurement, and operations, while leveraging off-peak electricity tariffs.}

\textcolor{blue}{While robust optimization is computationally efficient and well-suited to data-scarce environments, its reliance on worst-case scenarios often leads to overly conservative solutions that fail to capture the probabilistic nature of real-world uncertainty \citep{birge2011introduction}. To address these shortcomings, distributionally robust optimization (DRO) incorporates distributional information to balance robustness with performance, offering a more flexible approach to decision-making under uncertainty.}

\textcolor{blue}{DRO extends beyond the limitations of both stochastic programming and robust optimization by considering the worst-case probability distribution within an ambiguity set defined by available data \citep{delage2010distributionally}. Unlike stochastic programming, which relies on predefined probability distributions, and robust optimization, which assumes worst-case parameter realizations, DRO reduces conservatism while maintaining computational tractability. This makes it particularly effective for problems with limited or uncertain data \citep{mohajerin2018data,delage2010distributionally}.}

\textcolor{blue}{A notable extension of DRO is the distributionally robust chance constraint (DRCC) approach, which applies the DRO framework to probabilistic constraints. DRCC ensures that constraints are satisfied with a specified probability across all distributions in the ambiguity set, effectively balancing reliability and flexibility. This feature is particularly advantageous in addressing the challenges of sparse and uncertain data, where traditional methods may either fail to generalize or result in overly conservative designs. In this study, we adopt the DRCC approach to incorporate empirical observations into the optimization process, enabling a robust and adaptive framework for handling energy consumption variability in charging infrastructure design.}

\textcolor{blue}{\Cref{LR-summary}  presents the most relevant studies focusing on the incorporation of uncertainty into the optimization of BEB networks. While robust optimization methods with budget parameters to control conservatism are widely used, they can be overly conservative. To address this, our study introduces a novel application of the data-driven DRO approach for designing BEB charging infrastructure. Furthermore, this study evaluates how these techniques impact system performance, focusing specifically on cost efficiency and the feasibility of the resulting solutions.}

\begin{table}[htbp]
    \centering
    \renewcommand{\arraystretch}{1.5}
    \captionsetup{justification=raggedright,singlelinecheck=false,skip=2pt}
    \caption{\textcolor{blue}{Relevant studies in BEB planning design under uncertainty}}
    \tiny
    \textcolor{blue}{
    \begin{tabular}{p{4cm}p{0.75cm}p{1cm}p{1.15cm}p{0.75cm}p{1.5cm}p{1.75cm}}
    \toprule
       \multirow{2}{*}{\textbf{Authors}} & \multicolumn{2}{l}{\textbf{Decision level}}&\textbf{Data feature}& \multicolumn{3}{l}{\textbf{Modeling approach}}\\
 \cline{5-7}& & && SP& \multicolumn{2}{l}{Robust optimization}\\
       \cline{2-3} \cline{6-7}&Strategic&Tactical/ operational&& & Uncertainty set with budget&Data-driven DRO\\
       \midrule
       \cite{liu2018-TRR-planning} & \checkmark & \checkmark &Range&&\checkmark&\\
\cite{Kun-an2020battery} & \checkmark& &Scenario& \checkmark& &\\
       \cite{iliopoulou2021robust-BoU-partc} & \checkmark& &Range&&  \checkmark &\\
       \cite{bai2022robust} & \checkmark&\checkmark &Range& &  \checkmark  &\\
       \cite{hu2022joint} & \checkmark& \checkmark&Range& &  \checkmark  &\\
       \cite{huang2023robust} & &\checkmark&Range& & \checkmark  &\\
       \cite{cited-shadi_zhou2023electric} & \checkmark& &Scenario& \checkmark&  &\\
       \cite{gairola2023optimization} & \checkmark& \checkmark&Range& &  \checkmark  &\\
       \cite{avishan2023electric} & &\checkmark  &Range& &  \checkmark  &\\
       \midrule
       \textbf{This paper} & \checkmark& &Range and Sparse data & &\checkmark   &\checkmark   \\
       \bottomrule
    \end{tabular}}\\
    Legend-- SP: Stochastic Programming; DRO: Distributionally Robust Optimization \label{LR-summary}
\end{table}

\section{Model Formulation}
\label{sec:MF}
\textcolor{blue}{We formulate the CID problem for the electrification of an existing bus network, addressing the challenges posed by sparse energy consumption data. The problem involves three core decisions: (1) selecting the type of charging station (e.g., fast-feeding or standard), (2) identifying optimal charging station locations, and (3) determining onboard battery capacities.}

\textcolor{blue}{This section begins by detailing the fundamental assumptions and components of the electrification problem. We then present the formulation of three models: the deterministic CID model, the BoU-CID, and the DRCC-CID. Each model adopts a unique approach to handling energy consumption data. In \Cref{sec:NominalModal}, the (deterministic) CID model based on expected energy consumption is presented.  Next, \Cref{sec:BoUmodel} discusses the BoU-CID model which incorporates a range of uncertain values. Finally,  \Cref{sec:drccmodel}, presents the DRCC-CID model, which directly integrates empirical observations into a robust optimization framework.
}

\textbf{Cost structure.} The system cost structure includes infrastructure costs for charging stations and onboard batteries. Energy consumption uncertainty affects charge-discharge patterns, impacting battery longevity and system costs in long term. Therefore, we perform analysis on battery lifetime under different approaches.

\textbf{Network composition.} The bus network consists of several lines, denoted by the set $k \in K$. Each line $k$ requires a known number of buses. \textcolor{blue}{The BEB fleet size of bus line $k$ is denoted by $\gamma_k$.} All buses on line $k$ start and end their trips at a terminal, denoted by $o_k$, which is equipped with charging stations. Bus lines can be either bidirectional or one-way, with terminals at both ends. Each line is represented by an ordered set of stops, $S_k$, and the set of all stops is $S = \cup_{k \in K} S_k$.

\textbf{Charging infrastructure and battery specifications.} As mentioned before, we consider two types of charging station: fast-feeding and standard stations, denoted by $t\in T$. \textcolor{blue}{Fast-feeding stations provide rapid charging, enabling batteries to recharge up to 80\% in less than 20 seconds, making them ideal for intermediate stops where minimizing operational delays is critical. These systems are engineered to meet stringent safety standards, allowing safe operation even with passengers on board \citep{auge2015keynote}. Standard charging stations, in contrast, operate at a slower rate, with charging levels determined by the dwelling time and charger power. While these stations are better suited for end-of-line or depot locations, their slower charging rate can pose challenges for high-frequency operations.
In this study, we assume that charging stations can be installed at any bus stop. However, logistical constraints, such as limited grid capacity or local policy restrictions, may limit the set of feasible locations in practice. These constraints can be incorporated into the model by excluding specific stops from the set of potential charging locations}

Buses operating on the same line are assumed to have the same battery capacity. They must maintain a minimum energy level $\underline{b}$ and a maximum energy level $\bar{b}$ at all times. Additionally, BEBs are assumed to leave the terminal with maximum energy levels. 

\textbf{Energy consumption.} We assume that energy consumption between consecutive stops is a random variable, with some observations available from field experiments. \textcolor{blue}{To address energy consumption uncertainty, we do not constrain the analysis to specific factors such as travel time or BEB speed. Instead, the model is designed to capture variability arising from unobserved factors, such as environmental conditions or operational inconsistencies, that influence energy consumption during a route. This flexibility ensures that the model remains applicable across diverse scenarios where detailed information on uncertainty drivers may be unavailable. Building on observed energy consumption data, we propose a set of robust and data-driven modeling approaches to systematically incorporate uncertainty into the charging infrastructure design process.}

\subsection{Formulation Using Mean Energy Consumption Estimates}
\label{sec:NominalModal}
This model utilizes expected values to estimate uncertain energy consumption, serving as a baseline to gauge operational efficiencies. \textcolor{blue}{The notations used in the modeling formulation are summarized in \Cref{tab:notation}}. The objective function \textcolor{blue}{of all the models} (\ref{obj}) aims to minimize the total cost, compromising the installation cost of \textcolor{blue}{type $t \in T$ of charging stations (denoted by $\alpha_{t}$)} and onboard battery costs per kWh (denoted by $\beta$). \textcolor{blue}{The BEB fleet size for each line $k\in K$ is given and denoted by $\gamma_k$}. The binary variable $x_{st}$ is set to 1 if a charging station type $t \in T$ is installed at stop $s$. Since different bus lines can share a single charging station, the index $k\in K$ is omitted for this variable. The battery capacity for each bus line is represented by the continuous variable $z_k$. 
\begin{align}
    \min_{x,z} \quad \sum_{s \in S} \sum_{t \in T} \alpha_{t} x_{st} + \sum_{k \in K} \beta \gamma_k z_{k} \label{obj}
\end{align}

\begin{table}[htbp]
    \centering
    \renewcommand{\arraystretch}{1.5} 
    \captionsetup{justification=raggedright,singlelinecheck=false,skip=2pt} 
    \label{tab:notation}
    \caption{\textcolor{blue}{List of notations used in modeling formulation}}
    \tiny
    \textcolor{blue}{
    \begin{tabular}{p{1cm}p{12cm}}
    \toprule
        \textbf{Notation} & \textbf{Description}\\
        \midrule 
         \textbf{Sets}&\\
         \midrule
          $K$&Set of bus lines, indexed by k\\
 $S$&Set of all stops. $S_k$ denotes the stops on bus line $k$. The terminal (depot) of bus line $k$ is indexed as $o_k$\\
 $T$&Type of charging stations, indexed by t\\
 $N$&Sample size of observed energy consumption data, indexed by j\\
 \hline
 \textbf{Parameters}&\\
 \hline
 $\alpha_t$&Installation cost of charging station type $t$ (\euro)\\
 $\beta$&Unit cost per kWh of BEB battery capacity (\euro per kWh)\\
 $\gamma_k$&The BEB fleet size for bus line $k$\\
 $P_t$&Power capacity of charging station type $t$\\
 $\Delta_{ks}$&Dwell time of the BEB at stop $s$ on bus line $k$\\
 $\bar{b}$, \underline{$b$}&Upper and lower Bounds of allowable battery usage (\%)\\
 $\bar{\mu}^k_s, \hat{\mu}^k_s$&Average energy consumption and maximum deviation from the average between stops $s-1$ and $s$ on bus line $k$\\
 $\mu^{kj}_s$&$j$-th observed energy consumption sample between stops $s-1$ and $s$ on bus line $k$\\
 $\Gamma^k_s$&Budget parameter controlling the conservatism level in the \textbf{BoU-CID} model\\
 $\theta^k_s$&Budget parameter controlling the conservatism level in the \textbf{DRCC-CID} model\\
 $\epsilon$&Risk parameter for the chance constraint in the \textbf{DRCC-CID}\\
 \hline
 \multicolumn{2}{l}{\textbf{Decision variables}}\\
 \hline
 $e_{ks}$&BEB energy level upon departure from stop $s$ on bus line $k$\\
 $x_{st}$&Binary variable equal to 1 if a charging station of type $t$ is installed at stop $s$, and 0 otherwise.\\
           $z_k$&Continuous variable representing BEB battery capacity (kWh)\\
        \bottomrule
        \end{tabular}
        }
        \end{table}

Constraints for charging stations, energy flow, and  battery energy levels are presented below. In the remainder of this section, we explain each family of constraints separately.

\textbf{Charging station.}
In this problem, we assume that at most one type of charging station can be installed at each stop. This condition is enforced by Constraint \ref{one_type_charging}.
\begin{align}
    &\sum_{t \in T} x_{st} \leq 1, & \forall s \in S \label{one_type_charging}\\
    &x_{st} \in \{0,1\}, & \forall s \in S, \forall t \in T \label{binary_var}  
\end{align}

\textbf{Energy flow.} The variable $e_{ks}$ represents the energy level of the BEB on line $k$ as it departs stop $s$. Constraint (\ref{change_in_battery_intermediate}) governs the energy flow between stops $s-1$ and $s$. It ensures that the energy level at departure from stop $s$ is the energy level at departure from stop $s-1$ minus the expected energy consumption between these stops, \textcolor{blue}{denoted by $\bar{\mu}^k_s$, plus} any energy gained at stop $s$. The power of the charging station type $t$ is denoted by $P_t$, and $\Delta_{ks}$ indicates the dwelling time at stop $s$ for bus line $k$. \textcolor{blue}{For fast-feeding charging stations, the BEB departs with a fully charged battery regardless of the dwelling time, which is modeled as $P_{FF} = \bar{b} z_k$ .} BEBs must depart the terminal with maximum energy levels, given by $e_{ko_k} = \bar{b}z_k, \quad \forall k \in K$. 

\textcolor{blue}{Using the energy flow Constraint (\ref{change_in_battery_intermediate}), $e_{ks}$ is computed cumulatively from the terminal, as defined in Constraint (\ref{eq:cumu_e_leav}). This constraint represents the BEB energy level upon departure from stop $s$ on line $k$. It is expressed as the maximum battery level at the terminal, $\bar{b}z_k$, minus the cumulative average energy consumption across all preceding stops, $\sum_{i=o_k}^{s} \bar{\mu}^k_i$, plus the cumulative energy gained from charging at those stops, $\sum_{t \in T}\sum_{i=o_k}^{s} P_t \Delta_{ki} x_{it}$.}
\begin{align}
& e_{ks} \leq e_{ks-1} - \bar{\mu}^k_s + \sum_{t \in T} P_t \Delta_{ks} x_{st}, & \forall k \in K, s \in S_k \label{change_in_battery_intermediate}\\
&   e_{ks} = \bar{b} z_{k} -\sum_{i=o_k}^{s} \bar{\mu}^k_i + \sum_{t \in T}\sum_{i=o_k}^{s} P_t \Delta_{ki} x_{it},   &\forall k \in K, s\in S_k \label{eq:cumu_e_leav}
\end{align} 

\textbf{Battery energy level.} The BEB energy level at all stops must not exceed the maximum allowed, expressed as $e_{ks} \leq \bar{b} z_{k}, \quad \forall k \in K, s \in S_k\setminus o_k$. \textcolor{blue}{When applied to Constraint \ref{eq:cumu_e_leav}, the  term $\bar{b} z_{k}$ cancels out, simplifying the constraint. This ensures that the cumulative energy gained through charging up to stop $s$ does not exceed the cumulative average energy consumption up to stop $s$, thereby preventing overcharging and adhering to the maximum allowed battery capacity, as shown in Constraint (\ref{eq:uncertain1}).} Constraint (\ref{eq:uncertain2}) ensures that the energy level \textcolor{blue}{upon arrival} at any stop $s$ remains above the minimum required. The non-negativity of the battery capacity is shown in Constraint (\ref{non_negativity_cap}). 
\begin{align}
  & \sum_{t \in T} \sum_{i=o_k}^{s} P_t \Delta_{ki} x_{it} \leq \sum_{i=o_k}^{s} \bar{\mu}^k_{i} &\forall k \in K, s\in S_k \label{eq:uncertain1}\\
& \bar{b} z_{k}  -\sum_{i=o_k}^{s} \bar{\mu}^k_{i} + \sum_{t \in T} \sum_{i=o_k}^{s-1} P_t \Delta_{ki} x_{it} \geq \underline{b} \ z_{k}  &\forall k \in K, s\in S_k \label{eq:uncertain2}\\
& z_{k} \in \mathbb{R}^+, & \forall k \in K \label{non_negativity_cap}
\end{align}
The \textcolor{blue}{(deterministic) CID} is formulated as follows.
\begin{align}
\bm{CID} =  &\min_{x,z} \quad \sum_{s \in S} \sum_{t \in T} \alpha_{t} x_{st} + \sum_{k \in K} \beta \gamma_k z_{k} \notag&\\
\notag
&\text{subject to: } \ (\ref{one_type_charging},\ \ref{binary_var})\notag&\\ 
& \qquad \qquad \quad \ (\ref{eq:uncertain1}-\ref{non_negativity_cap})\notag&
\end{align}

\textcolor{blue}{CID} provides electrification solution by using the expected value of energy consumption. However, various contextual factors, such as weather, driving style, and passenger loads, cause deviations from this average consumption. If the true distribution of energy consumption is skewed, relying solely on average estimates may misrepresent the actual conditions \textcolor{blue}{and result in a charging infrastructure design that fails to provide feasible and reliable trips. This highlights the need for more comprehensive formulations that incorporate energy consumption uncertainty into the BEB CID. Consequently, two robust models are introduced and discussed below.}

\subsection{Formulation based on uncertainty range}
\label{sec:BoUmodel}
In this approach, we consider energy consumption uncertainty as a range. This strategy is applicable when the range of energy consumption data is known to decision-makers. For two subsequent stops $s-1$ and $s$, \textcolor{blue}{the uncertainty box is defined by the expected energy consumption ($\bar{\mu}^k_s$) and its maximum deviation ($\hat{\mu}^k_s$). The deviation can be expressed as $\hat{\mu}^k_s= \bar{\mu}^k_{s} + \omega\bar{\mu}^k_{s}$. Therefore, the energy consumption belongs to the range of $[\bar{\mu}^k_{s}-\omega\bar{\mu}^k_{s},\bar{\mu}^k_{s}+\omega\bar{\mu}^k_{s}]$, where $\omega \in[0,1]$ indicating possible deviations up to the average energy consumption value \cite{bertsimas2003robust}}. Negative deviations are disregarded as they do not impact the bus system's service level. Consequently, the uncertainty range is defined as:   $[\bar{\mu}^k_{s},\bar{\mu}^k_{s}+\omega\bar{\mu}^k_{s}]$.

We standardized the \textcolor{blue}{random energy consumption parameter ($\mu^k_s$)} by defining $\upvarphi^k_{s} = \frac{\mu^k_{s} - \bar{\mu}^k_{s}}{\hat{\mu}^k_{s}}$. Without constraints on $\upvarphi^k_{s}$, this approach can be overly conservative, allowing maximal deviations in energy consumption between all stops on the bus line. To mitigate this, we impose $\Gamma^k_{s}$ as a budget to set an upper limit for energy consumption deviations up to stop $s$ on bus line $k$. We define the box uncertainty set for $\upvarphi^k_{n}$ as $\mathcal{C}(\upvarphi^k_n)$: 

\begin{align*}
    \mathcal{C}(\upvarphi^k_n) = \{\upvarphi^k_i \in[0,1] | \ \sum_{i=o_k}^s\upvarphi^k_{i}\leq \Gamma^k_{s}, \quad \forall k\in K\}
\end{align*}
The  parameter $\Gamma^k_{s}$ belongs to the interval $[0,|s|]$. This parameter indicates the cumulative allowed energy consumption deviation up to stop $s$. A value of $\Gamma^k_{s}= |s|$ means that the energy consumption is at its maximum level. The value of $\Gamma^k_{s}$ depends on decision-maker opinion and indicates the robustness \textcolor{blue}{and level of conservatism} of the solution. \textcolor{blue}{We will now consider the range of the energy consumption random parameter in (\ref{eq:uncertain1}) and (\ref{eq:uncertain2}) using the defined box uncertainty set $\mathcal{C}(\upvarphi^k_n)$ as shown below}:
\begin{align}
    & \sum_{t \in T} \sum_{i=o_k}^{s} P_t \Delta_{ki} x_{it} \leq \sum_{i=o_k}^{s} \bar{\mu}^k_{i} +\sum_{i=o_k}^{s} \hat{\mu}^k_{i} \upvarphi^k_{i},  &\forall k \in K, s\in S_k  \label{eq:1stBoU}\\
   & \bar{b} z_{k}  -\sum_{i= o_k}^{s} \bar{\mu}^k_{i} -\sum_{i= o_k}^{s} \hat{\mu}^k_{i} \upvarphi^k_{i} + \sum_{t \in T} \sum_{i= o_k}^{s-1} P_t \Delta_{ki} x_{it} \geq \underline{b} \ z_{k}  &\forall k \in K, s\in S_k \label{eq:2ndBoU}
\end{align}
\textcolor{blue}{The above constraints provide a formulation where the energy consumption is estimated using a box uncertainty set. This strategy aims to optimize BEB CID to be robust against maximum energy consumption deviation up to stop $s$, ($\sum_{i= o_k}^{s} \hat{\mu}^k_{i}\upvarphi^k_{i}$), within a predefined budget from the terminal up to stop s as $\Gamma^k_{s}$.} This term in the above constraints is represented by the following optimization model:
\begin{align}
\max_\upvarphi & \sum_{i= o_k}^{s} \hat{\mu}^k_{i}\upvarphi^k_{i} \notag \\
\text{s.t.} & \sum_{i= o_k}^{s} \upvarphi^k_{i} \leq \Gamma^k_{s}, \quad  \forall s\in S_k \notag \\
& 0 \leq \upvarphi^k_{s} \leq 1, \quad \forall s\in S_k
\label{primal}
\end{align}
The dual formulation can be written as:
\begin{align}
\min_{u,v} \quad & \Gamma^k_{s}u^k + \sum_{i= o_k}^{s} v_i^{k} \notag \\
\text{s.t.}\quad &   u^k + v_s^{k} \geq \hat{\mu}^k_{s} , \quad \forall s\in S_k \notag \\
& u^k, v_s^{k} \geq 0 , \qquad\forall s\in S_k
\end{align}

\textcolor{blue}{Consequently, Constraints (\ref{eq:1stBoU}) and (\ref{eq:2ndBoU}) are reformulated as (\ref{BoU-cumulative}) and (\ref{constr:BoU3}).}
\begin{align}
     & \sum_{t \in T} \sum_{i=o_k}^{s} P_t \Delta_{ki} x_{it} \leq \sum_{i=o_k}^{s} \bar{\mu}^k_{i} + \Gamma^k_{s}u^k + \sum_{i=o_k}^s v_i^{k}, \qquad \qquad \forall k \in K, s\in S_k  \label{BoU-cumulative}&\\
   & (\bar{b}-\underline{b}) z_{k} + \sum_{t \in T} \sum_{i=o_k}^{s-1} P_t \Delta_{ki} x_{it} \geq \sum_{i=o_k}^{s} \bar{\mu}^k_{i} +\Gamma^k_{s}u^k + \sum_{i=o_k}^s v_i^{k}, \quad \forall k \in K, s\in S_k  \label{constr:main-BoU}&\\
   & u^k + v_s^{k} \geq \hat{\mu}^k_{s}, \qquad \qquad \qquad \qquad \qquad  \qquad\qquad \qquad \forall k \in K,s \in S_k \label{constr:BoU2}&\\
    & u^k, v_s^{k} \geq 0, \qquad \qquad \qquad \qquad \qquad \qquad  \qquad \qquad \qquad \forall s\in S_k&
    \label{constr:BoU3}
\end{align} 

\textcolor{blue}{The deterministic equivalent counterpart of the box uncertainty robust model, i.e., BoU-CID,  is as follows.}
\begin{align}
\bm{BoU-CID} =  &\min_{x,z} \quad \sum_{s \in S} \sum_{t \in T} \alpha_{t} x_{st} + \sum_{k \in K} \beta \gamma_k z_{k} \notag&\\
\notag
&\text{subject to: } \ (\ref{one_type_charging},\ \ref{non_negativity_cap},\ \ref{binary_var})\notag&\\ 
& \qquad \qquad \quad \ (\ref{BoU-cumulative}-\ref{constr:BoU3})\notag&
\end{align}

\subsection{Empirical Data-Driven Uncertainty Modeling}
\label{sec:drccmodel}
\textcolor{blue}{In the previous section,} we demonstrated the importance of accounting for values above the average consumption to ensure robustness of the design and developed a tractable formulation based on this concept. However, as discussed in \cite{baron2011facility-robust,liu2019-EMS-DRCC}, this model may be vulnerable to inaccuracies in estimating the worst-case scenarios. An exceptionally high but less probable worst-case realization of energy consumption at a stop can disproportionately influence the design of the entire network. 

In this section, we introduce the third approach to handle the uncertainty in energy consumption, assuming access to a limited number of observational data points. The main idea is to construct a distributional set of energy consumption that integrates more data points beyond mere range information. By utilizing these observations to approximate the true probability distribution of energy consumption, this approach provides a solution robust against the worst probable distribution within a predefined budget. Consequently, it reduces reliance on extreme values and seeks to lower costs from  overly conservative designs.

\textcolor{blue}{We consider $N$ samples, each representing observations of energy consumption between consecutive bus stops on a specific bus line.} The $j$-th energy consumption sample between stop $s-1$ and $s$ on bus line $k$ is noted by $\mu^{kj}_s$, where $\{j= 1,...,N\}$. The goal of this approach is to find a robust design that minimizes the electrification costs utilizing $N$ energy consumption samples. \textcolor{blue}{It ensures that the cumulative energy consumption up to stop $s$ on line $k$ falls within a decision-dependent safety set with high probability ($1-\epsilon$) across all potential distribution ($\mathds{P}$), as shown in Constraints (\ref{DRCC-constr_cumulative}) and (\ref{const2-DRCC}). In these constraints, the average energy consumption in \textbf{CID} is replaced by the observed values of energy consumption. The related risk parameter is shown by $\epsilon$.}
\begin{align}
    & \mathds{P} \Biggl[ \sum_{t \in T} \sum_{i= o_k}^{s} P_t \Delta_{ki} x_{it} \leq \sum_{i= o_k}^{s} \mu^{kj}_{i}, \quad \forall k \in K, s\in S_k,j\in\{1,...,N\}  \Biggr] \geq 1 - \epsilon, \notag\\
    &\qquad \qquad \qquad \qquad \qquad \qquad\qquad \qquad \qquad\qquad \qquad \qquad \forall \mathds{P}\in \mathcal{F}(\theta^k_s) \label{DRCC-constr_cumulative}&\\
    & \mathds{P} \Biggl[ (\bar{b}-\underline{b}) z_{k} +\sum_{t \in T} \sum_{i= o_k}^{s-1} P_t \Delta_{ki} x_{it} \geq \sum_{i= o_k}^{s} \mu^{kj}_{i},  \forall k \in K, s\in S_k, j\in\{1,...,N\} \Biggr] \geq 1 - \epsilon,\notag\\
    &\qquad \qquad \qquad \qquad \qquad \qquad\qquad \qquad \qquad\qquad \qquad \qquad \quad\forall \mathds{P}\in \mathcal{F}(\theta^k_s)
    \label{const2-DRCC}
\end{align}
Given the $N$ observations of energy consumption data, accurately estimating their true distribution, which is influenced by multiple factors, is complex. \textcolor{blue}{To address this, we define a distributional uncertainty set, or ambiguity set  (denoted by\( \mathcal{F}(\theta^k_s) \)), that includes all possible distributions for energy consumption (\( \mathds{P} \)) within a specified distance (\( \theta^k_s \)) from a reference distribution, denoted by $\hat{\mathds{P}}$.} The reference distribution is modeled as a discrete uniform distribution calculated from the energy consumption observations, expressed as \( \hat{\mathds{P}} = \frac{1}{N} \sum_{j=1}^N \delta_{\mu^{kj}_s} \), where $\delta$ represents the Dirac delta function for each energy consumption observation $\mu^{kj}_s$. \textcolor{blue}{The ambiguity set is defined by $\mathcal{F}(\theta^k_s) = \{ \mathds{P}: d_W{(\hat{\mathds{P}},\mathds{P})} \leq \theta^k_s \}$, with $d_W$ being the Wasserstein metric measuring distributional distance.} This metric effectively forms a ball of radius $\theta^k_s \geq 0$ centered on the empirical distribution $\hat{\mathds{P}}$.

The decision-dependent safety set, denoted as ${\mathcal{S}}(X)$, ensures that constraints (\ref{DRCC-constr_cumulative}) and (\ref{const2-DRCC}) are met. This  allows the model to probe the Boundary of the safety set and identify the worst probable distribution of energy consumption influenced by the parameter $\theta^k_s$. The parameter $\theta^k_s$ serves as a budget for displacement, dictating the extent to which observed energy data up to stop $s$ on bus line $k$ can shift toward the Boundary of its safety set. \textcolor{blue}{The decision maker determines the value of $\theta^k_s$, which sets the size of the ambiguity set. A larger $\theta^k_s$ value provides more budget to transport an energy consumption sample to its worst value, leading to more conservative solutions. If $\theta^k_s = 0$, the model considers only the average value of the energy consumption data, i.e., \( \bm{CID} \).}

We denote unsafe set by $\bar{\mathcal{S}}(X)$, representing the infeasible region of the model, which is the complement of ${\mathcal{S}}(X)$. \textcolor{blue}{To shift the observed data to the border of its safety set, the distance of the uncertain energy consumption data $\mu^{kj}_s$, to the unsafe set is defined and denoted by $\text{dist}(\mu^{kj}_s, \bar{\mathcal{S}}(X))$ as follows:}
\begin{align}
    \text{dist}(\mu^{kj}_s, \bar{\mathcal{S}}(X)) = \max\{0, &\min ( (\sum_{i=o_k}^{s} \mu^{kj}_{i}-\sum_{t \in T} \sum_{i=o_k}^{s} P_t \Delta_{ki} x_{it}) \ \textbf{,} \notag\\
    &((\bar{b}-\underline{b}) z_{k} +\sum_{t \in T} \sum_{i= o_k}^{s-1} P_t \Delta_{ki} x_{it} - \sum_{i= o_k}^{s} \mu^{kj}_{i}))) \} \notag\\
    &\qquad \qquad \qquad \qquad \qquad \qquad \forall j\in\{1,...,N\}\label{non-convex-dist-unsafe}
\end{align}

For the energy consumption data \( \mu^{kj}_s \) within the constraints (\ref{DRCC-constr_cumulative}) and (\ref{const2-DRCC}), the distance to the unsafe set is calculated as the minimum distance to any of these  constraints. This distance represents the required adjustment to move an observed energy consumption value towards its worst-case value. Energy consumption values that already fall outside the safety set are assigned a zero distance, as they are already in the unsafe set. The chance constraint ensures that the constraints are jointly feasible with a probability of at least $(1-\epsilon)$ for data points within the safety set. 

After calculating the distance of each observation $\mu^{kj}_s$ to the unsafe set, we sort these samples by increasing distance and reorder them using the index $\pi_j(X)$. The sorted energy consumption is then represented by $\mu^{k\pi_j(X)}_s$. This index arranges the observations based on their proximity to the unsafe set, contingent upon the decision variable $X$. The decision variable in the index represent the design with different robustness settings. The \textit{displacement budget}, $\theta^k_s$, is a key robustness setting that determines how many and which energy consumption data points are shifted towards the unsafe set to derive the worst-case energy consumption distribution as per constraint (\ref{control_dist}). The energy consumption data points are moved in the order specified by $\pi_j(X)$, within the limits of $\theta^k_s$, to approximate the worst probable distribution of energy consumption. Detailed proofs are provided in \cite{kuhn-chen2022data}. 

\begin{align}
    \sum _{j=1}^{\epsilon N}\text{dist}(\mu^{k  \pi_j(X)}_s, \bar{\mathcal{S}}(X)) \geq \theta^k_s N\label{control_dist}
\end{align}
Constraint (\ref{control_dist}) ensures that the distances of observations from the unsafe set are adjusted according to the robustness parameter, $\theta^k_s$. Specifically, this equation selects the first $\epsilon N$ energy consumption samples nearest to the unsafe set and moves them towards the Boundary, determining the worst energy consumption distribution. 

The calculation of distances to the unsafe set involves Constraints (\ref{non-convex-dist-unsafe}) and (\ref{control_dist}). To eliminate the partial sum operator and the permutation index $\pi_j(X)$, Constraint (\ref{control_dist}) is replaced by equations (\ref{eps}) and (\ref{linear_dist}), without approximation, but at the cost of introducing additional variables $q^k \geq 0$ and $r^{kj}_s \geq 0$. These variables are  dual variables.
\begin{align}
    & \epsilon N q^{k} - \sum_{i=o_k}^s \sum_{j\in\{1,...,N\}} r^{kj}_i \geq \theta^k_s N  \qquad \qquad \quad \forall k\in K, s\in S_k \label{eps}\\
    &\text{dist}(\mu^{kj}_s, \bar{\mathcal{S}}(X)) \geq q^{k} - \sum_{i= o_k}^n r^{kj}_i \qquad \qquad  \forall k\in K,, s\in S_k, j\in\{1,...,N\} \label{linear_dist}   
\end{align}

\textcolor{blue}{The final step in the DRCC-CID formulation is the linearization of the distance function (\ref{linear_dist}).}Following \cite{kuhn-chen2022data}, we introduce a binary variable $y^{kj}_s$, where $y^{kj}_s= 1$ indicates that the sample $\mu^{kj}_{s}$ does not satisfy the chance constraint, and vice versa. 
\begin{align}
    & M(1- \sum_{i= o_k}^s y^{kj}_i )\geq q^{k} - \sum_{i= o_k}^s r^{kj}_i \qquad \quad \forall k\in K,s\in S_k, j\in\{1,...,N\} \label{upper_BoUnd_yj}\\
    & - \sum_{t \in T} \sum_{i= o_k}^{s} P_t \Delta_{ki}\ x_{it} + \sum_{i= o_k}^{s} \mu^{kj}_{i} + M\ \sum_{i= o_k}^s y^{kj}_i \geq q^{k} - \sum_{i= o_k}^s r^{kj}_i \notag\\
    &\qquad \qquad \qquad \qquad \qquad \qquad\qquad \qquad\forall  k\in K, s \in S_k, j\in\{1,...,N\} \label{DRCC-new_constr5}\\
    &(\bar{b}-\underline{b}) z_k + \sum_{t\in T}\sum_{i= o_k}^{s-1} P_t \Delta_{ki}x_{it}- \sum_{i= o_k}^s \mu^{kj}_{i} + M\ \sum_{i= o_k}^s y^{kj}_i \geq q^{k} - \sum_{i= o_k}^s r^{kj}_i \notag\\
    & \qquad \qquad \qquad \qquad \qquad \qquad \qquad \qquad \forall  k\in K, s \in S_k, j\in\{1,...,N\} \label{drjcc-const9}\\
    & y^{kj}_s \in \{0,1\}^N, \ r^{kj}_s,\ q^{k}\geq0, \qquad \qquad \qquad \qquad \quad \quad \forall k\in K, s\in S_k\label{eq-26}
\end{align}
Constraints (\ref{eps}) and 
(\ref{upper_BoUnd_yj}--\ref{eq-26}) limit the distance of observations to the unsafe set, ensuring that the probability of not satisfying the chance constraint is smaller than $\epsilon$. Constraints (\ref{upper_BoUnd_yj}--\ref{drjcc-const9}) include a sufficiently large constant $M \in \mathbb{R}^+$, ensuring that the value of $\text{dist}(\mu^{kj}_s, \bar{\mathcal{S}}(X))$ is always non-negative.

The DRCC Constraints (\ref{DRCC-constr_cumulative}) and (\ref{const2-DRCC}) can now be replaced by their equivalents (\ref{eps}) and (\ref{upper_BoUnd_yj}--\ref{eq-26}). \textcolor{blue}{The $\bm{DRCC-CID}$ model, is ultimately reformulated as an MILP and presented below:}

\begin{align}
\bm{DRCC-CID} = &\min_{y,r,q,x,z}  \sum_{s \in S} \sum_{t \in T} \alpha_{t} x_{st} + \sum_{k \in K} \beta \gamma_k z_{k} \notag\\
\notag
&\text{subject to:} \ (\ref{one_type_charging},\ \ref{binary_var},\ \ref{non_negativity_cap}),\notag&\\
    & \qquad \qquad \quad (\ref{eps}) \  \text{and} \ (\ref{upper_BoUnd_yj}-\ref{eq-26}) \notag
\end{align}

\section{Computational results}
\label{sec:results}
\textcolor{blue}{In this section, we outline the data structure required to implement the proposed models (Section \ref{sec:data-structure}) and evaluate their applicability and scalability. The performance of the CID, BoU-CID, and DRCC-CID models is first tested on a series of hypothetical grid networks (Section \ref{sec:hypo-grid}) to assess scalability and subsequently applied to the Rotterdam bus network (Section \ref{sec:Rdam-case}) to demonstrate their practical effectiveness. To compare the feasibility of the design solutions generated by the CID, BoU-CID, and DRCC-CID approaches, a Monte Carlo simulation procedure is detailed in Section \ref{sec:simulation-monte}. The impact of sample size on the performance of the DRCC-CID model is analyzed in Section \ref{sec:quantity-drcc}, while Section \ref{sec:sensitivity} presents sensitivity analyses on economic factors, including charging station costs and BEB fleet size.
The models are implemented using Python 3.9.13, with the Gurobi 9.5.2 optimizer serving as the MIP solver. All computational experiments were conducted on an Apple M1 with 16 GB of RAM, ensuring reliable and reproducible results.} 

\subsection{Data Structure}
\label{sec:data-structure}
\textcolor{blue}{This section details the specifications of the datasets used in our study, including both case study and synthetic data. These data sets are used to evaluate the effectiveness of the proposed approach.}

\begin{figure}[htbp]
    \centering
    \begin{subfigure}[b]{\textwidth} 
        \includegraphics[width=\textwidth]{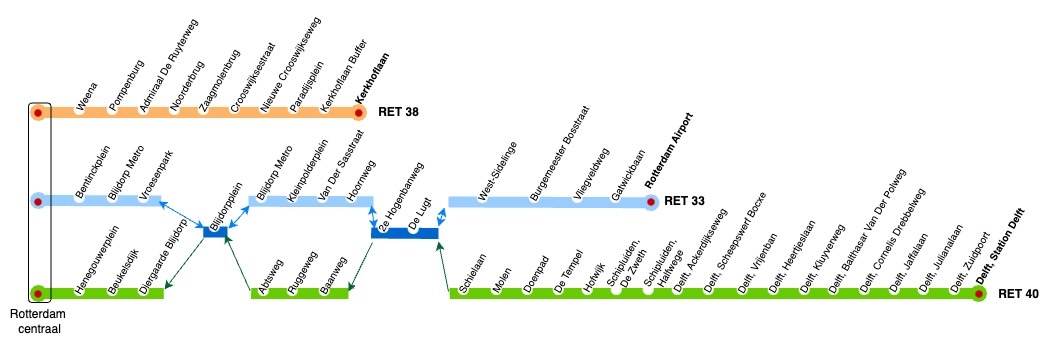}
        \caption{Graphical depiction of Rotterdam bus network for bus lines 33, 38, and 40}
        \label{fig:bus-network}
    \end{subfigure}
    \begin{subfigure}[b]{0.3\textwidth}
        \centering
        \includegraphics[width=\textwidth]{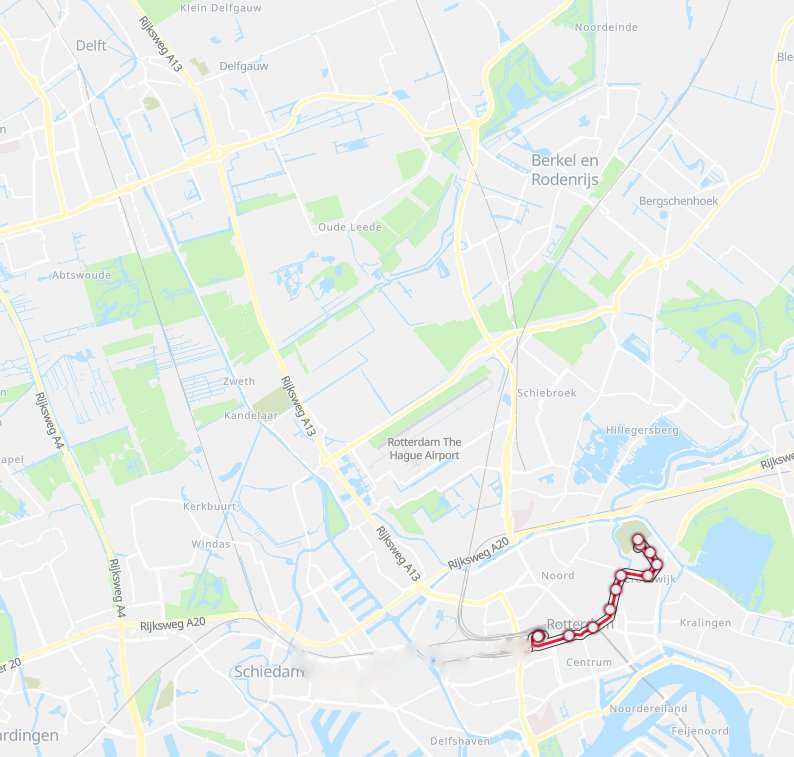}
        \caption{Part of bus line 38}
        \label{fig:busline38}
    \end{subfigure}
    \hfill
    \begin{subfigure}[b]{0.3\textwidth}
        \centering
        \includegraphics[width=\textwidth]{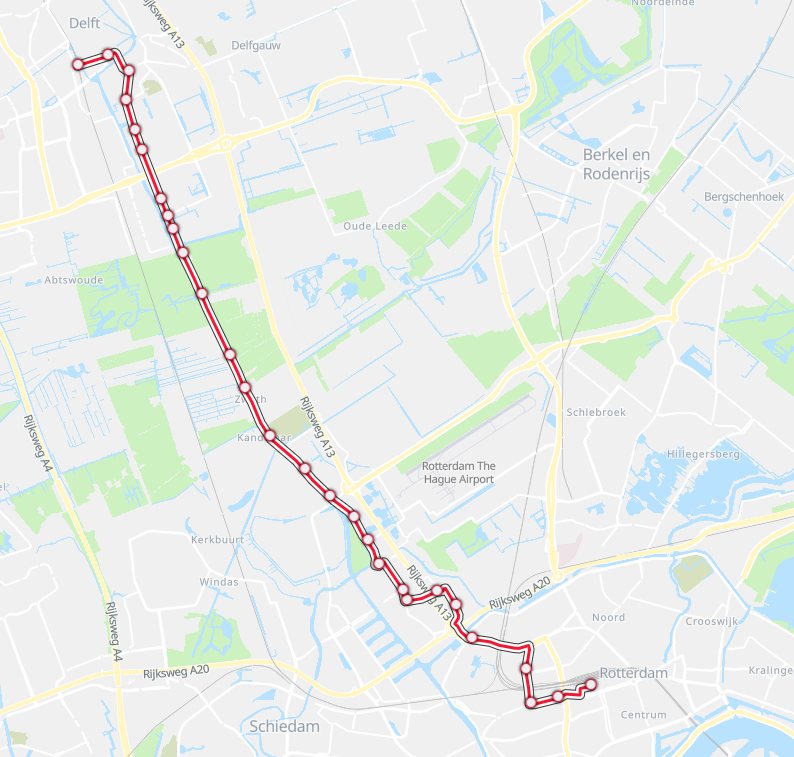}
        \caption{Bus line 40}
        \label{fig:busline40}
    \end{subfigure}
    \hfill
    \begin{subfigure}[b]{0.3\textwidth}
        \centering
        \includegraphics[width=\textwidth]{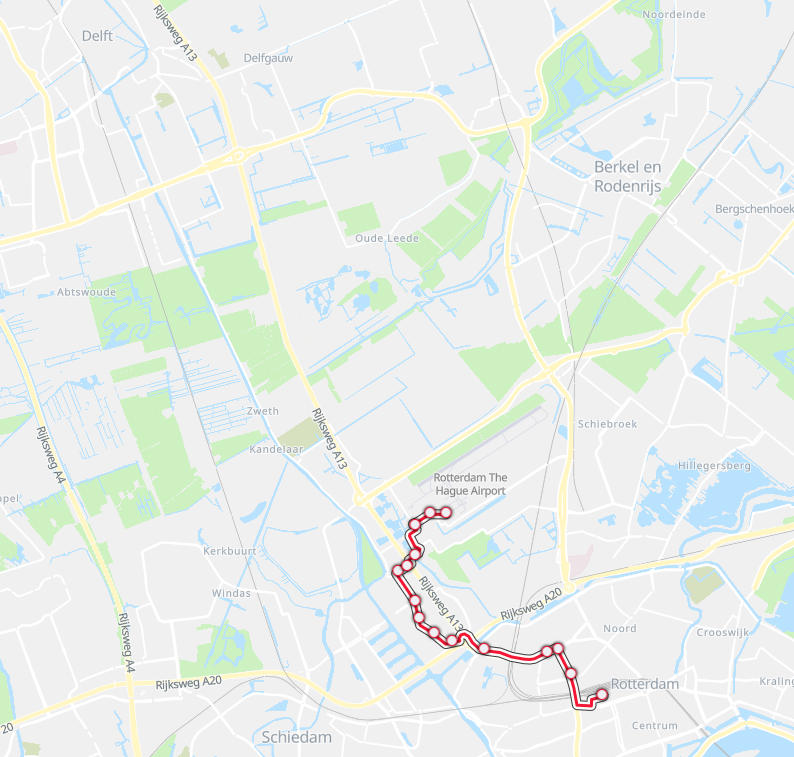}
        \caption{Bus line 33}
        \label{fig:busline33}
    \end{subfigure}
    \caption{Rotterdam bus lines depicted on the map}
    \label{fig:map-bus-network}
\end{figure}

\textcolor{blue}{\textbf{Hypothetical grid network structure.} To test the scalability of the proposed models, a $10\times 10$ grid network is generated to represent various bus network structures. This grid comprises 100 nodes located at intersections, serving as both bus stops and potential charging station locations. We simulate different scenarios with varying numbers of bus lines ($k\in\{5, 25, 45\}$), and bus stops per line ($s_k\in \{5, 25, 45\}$). In each scenario, number of bus stops per line are predefined, and their locations are randomly generated on the grid, resulting in interconnected routes with varying inter-stop distances. All synthetic grid-based bus lines share common start and end depots. Energy consumption data is simulated based on the distance between stops, as detailed in a subsequent section. A schematic view of the simulated bus lines with different scales is provided in Appendix \ref{fig:hypo-grids}. The model's performance in finding exact optimal solutions is evaluated across these scenarios.}\\
\textcolor{blue}{\textbf{Selected bus lines in Rotterdam bus network.} In this case study, we select three distinct bus lines within the Rotterdam bus network to investigate the deployment of charging infrastructure for BEBs. Our objective is to examine how different modeling approaches impact bus lines with varying spatial characteristics. Specifically, we analyze bus lines 33, 38, and 40, which include both short-distance and intercity routes. All lines begin at the common terminal,  \textit{Rotterdam Centraal Station}, and may share stops with one another. Each line is serviced by a fleet of 10 buses. A schematic representation of the studied network is shown in Figure (\ref{fig:map-bus-network}).\\Bus line 38 is an intra-city route with 11 closely spaced stops. In contrast, Bus line 40 is an intercity route connecting Rotterdam and Delft, with 28 widely spaced stops. Bus line 33 serves a peripheral area with 15 stops. Both bus lines 33 and 38 are circular routes that start and end at \textit{Rotterdam Centraal Station}, while bus line 40 operates as a one-way route with distinct terminals at each end. Additionally, bus lines 33 and 40 share three stops; however, bus line 33 serves these stops in both directions, unlike bus line 40.}

\textcolor{blue}{\textbf{Charging station parameters and BEB specification.} The specifications of the two types of charging stations, fast-feeding and standard stations, and BEB batteries are summarized in \Cref{parameters-infrastructure}.  We assume that all buses within a given line have homogeneous battery capacities, although different bus lines may have varying capacities. In our model, buses are assumed to charge en route at bus stops, as in \citep{bai2022robust,auge2015keynote}.}

\textcolor{blue}{The electrification of bus networks involves two primary cost components: initial investment and long-term operational expenses related to battery life. After the design phase, we evaluate the long-term cost based on battery life, specifically the number of cycles until BEB battery failure, $N_{cycle}$, which depends heavily on the charge-discharge patterns defined by each modeling approach. Following \citep{zang2022battery_life}, the Depth of Discharge at stop $s$ on line $k$ ($\text{DOD}_s^k$) and the charged received at each stop $s$ ($\text{charge}_s^k$) are defined as follows:
\begin{align}
&\text{DOD}_s^k = \frac{(z_k - y_s^{'k})}{z_k},\  & y_s^{'k} \in [0, \bar{b} z_k]\\
&\text{charge}_s^k = \frac{z_k-Y_s^k}{z_k},\ & Y_s^k \in [0, \bar{b} z_k]  
\end{align} 
where $y_s^{'k}$ is the BEB remaining energy upon arrival at stop $s$ and $Y_s^k$ is the battery level upon departure from stop $s$. The difference $(Y_s^k - y_s^{'k} )$ indicates the charge received at stop $s$. Therefore, $N_{cycle}$ for each bus line $k$ can be deduced as follows, with constants $g_1 = 1331$ and $g_2 = 1.825$ for typical Li-ion batteries \cite{dallinger2013plug}.
\begin{align}
\label{eq:n-cycle}
N_{cycle}^k = g_1 \cdot \sum_s ((\text{DOD}_s^k)^{-g_2} + (\text{charge}_s^k)^{-g_2})
\end{align}
The cost-efficiency of each model is evaluated by computing the \textit{battery cost per cycle} for each bus line as:
\begin{align}
\label{eq:cost-per-cyc}
&\text{Cost per cycle} = \frac{\beta\times z_k}{N_{cycle}}
\end{align} }
\begin{table}[htbp]
    \centering
    \renewcommand{\arraystretch}{1.5} 
    \captionsetup{justification=raggedright,singlelinecheck=false,skip=2pt} 
    \caption{List of parameters used in the models}
    \tiny
    \begin{tabular}{llll}
    \toprule
       \textbf{Parameter description} & \textbf{Notation} & \textbf{Value} & \textbf{Reference}\\
       \midrule
        \textbf{Infrastructure related} \\
       \midrule
        Charging type & $T$ & \{\textit{Standard, Fast-feeding}\} & - \\
        \multirow{2}{*}{Cost of installing charging stations} & $\alpha_{SS}$ & 20,000 \euro & Assumption (sensitivity analysis) \\
                                                               & $\alpha_{FF}$ & 80,000 \euro & Assumption (sensitivity analysis) \\
        Unit cost of battery & $\beta$ & 1750 \euro/kWh & \cite{lajunen2014energy}, \cite{liu2018-TRR-planning} \\
        \multirow{2}{*}{Maximum Power of chargers} & $P_S$ & 100 kW & \cite{parameter-pihlatie2017requirements} \\
                                                   & $P_{FF}$ & 600 kW & \cite{parameter-pihlatie2017requirements} \\
        Lower percentage of battery usage limit & $\underline{b}$ & 20\% & \cite{shadi-azadeh2022electrification} \\
        Upper percentage of battery usage limit & $\bar{b}$ & 80\% & \cite{liu2018-TRR-planning} \\
        BEB dwell time at stops & $\Delta$ & N(20,4) Seconds & \cite{auge2015keynote}, \cite{shadi-azadeh2022electrification} \\
    \bottomrule
    \end{tabular}
    \label{parameters-infrastructure}
\end{table}

\textcolor{blue}{\textbf{Energy consumption synthetic data generation.} Various methods have been proposed for estimating energy consumption in BEBs, accounting for key determinants such as aerodynamic design, rolling resistance between tires and road surfaces, gravitational forces, and other influencing factors \citep{fiori2021microscopic,fontana2013optimal}. While these methods provide detailed insights into energy consumption fluctuations, they often rely on specific parameter estimation. By treating energy consumption as an exogenous input to the model, we reduce dependency on detailed parameter estimation. In this study, we address energy consumption uncertainty without focusing on specific underlying factors, the average energy consumption ($\bar{\mu}$) is simulated based on the distance between two consecutive stops, with a consumption rate of 1.3 kWh per km \citep{bai2022robust}. Google Maps is used to calculate the distances between each pair of stops for the Rotterdam bus lines. For the BoU-CID robust model, the maximum deviation from the mean energy consumption ($\hat{\mu}= \bar{\mu} + \omega\bar{\mu}$) is randomly generated where $\omega \in[0,1]$, similar to \citep{bertsimas2003robust}. This selection indicates that the energy consumption between these locations can deviate up to their average energy consumption value.} 

\textcolor{blue}{In the DRCC-CID the base scenario involves generating 100 energy consumption samples from a uniform distribution within the range $[\bar{\mu}, \hat{\mu}]$ for each pair of stops. Other parameters are set as follows: the risk tolerance $\epsilon$ is set to $0.1$ to satisfy the chance constraints, and $M$ is fixed at 25 after systematic evaluation.}

\textcolor{blue}{\textbf{Budget parameters for robust models.} It is important to note that the budget parameters in the BoU-CID and DRCC-CID robust models differ in their definitions. Although both parameters represent the conservatism level in the solution, their interpretations differ between the models. In the BoU-CID model, the budget parameter ($\Gamma$) represents the maximum allowable deviation from the average energy consumption value. In contrast, the budget parameter in the DRCC-CID model ($\theta$) serves as a displacement budget, governing the extent to which the observed energy data can shift toward the safety set Boundary. This approach constructs the worst-case distribution of energy consumption, rather than emphasizing the worst individual value.} 

\textcolor{blue}{In this study, two conservatism levels are selected for the robust models to represent low- and high- conservatism solutions, with $\Gamma = \{0.2, 0.8\}$ and $\theta = \{0.2, 0.8\}$. A sensitivity analysis examining the impact of other values for the budget parameters on the robust model solutions for the Rotterdam bus network, is presented Appendix \cref{tab:BoU-performance,tab:DRCC-performance}.}

\subsection{Hypothetical grid network}
\label{sec:hypo-grid}
\textcolor{blue}{To evaluate the performance of CID, BoU-CID, and DRCC-CID models on large-scale networks, the models are applied to the grid-based bus network structure described in \ref{sec:data-structure}. Input parameters, including charging station installation costs, battery capacity, and BEB fleet size, remain consistent with earlier definitions to assess model scalability. The charging infrastructure design problems are solved using the Gurobi optimizer.}

\textcolor{blue}{For the deterministic solution, which uses average energy consumption in the modeling, all tested network scenarios were solved to optimality in under 10 seconds. With energy consumption uncertainty included in the BoU-CID model, all scenarios were also solved to optimally. The runtime for the largest network (45 bus lines, each with 45 stops) increased to approximately 16 minutes. This highlights the efficiency of robust solutions in finding optimal solutions within reasonable computational time. }

\textcolor{blue}{ \Cref{tab:scalability-grid} summarizes the results for DRCC-CID across different network sizes and uncertainty levels.  For low conservatism  $(\theta=0.2)$, DRCC-CID found optimal solutions for all scenarios except for the network with 5 bus lines and 45 stops, indicating that scalability is more affected by the number of stops per line than by the number of lines. This is due to the increase in the number of data in dense bus lines.}
\begin{table}[h]
\centering
    \renewcommand{\arraystretch}{1.5} 
    \captionsetup{justification=raggedright,singlelinecheck=false,skip=2pt} 
    \centering
    \caption{\textcolor{blue}{Scalability test scenarios of the DRCC-CID model on grid-based bus networks}}
\tiny
\label{tab:scalability-grid}
\textcolor{blue}{
    \begin{tabular}{p{1.5cm}p{1cm}p{1cm}p{1.35cm}p{1cm}p{1.5cm}p{1.5cm}p{0.5cm}p{0.5cm}}
    \toprule
          \textbf{Uncertainty level} ($\theta$)&\multicolumn{2}{l}{\textbf{Grid-based bus network size}} &  \multicolumn{2}{l}{\textbf{Electrification costs} ($10^6$\euro)}&\textbf{Final upper bound}&\textbf{Final lower bound}& \textbf{Final gap} (\%)&\textbf{CPU time} (s)\\
          \cmidrule{2-5}
          &Number of bus lines& Number of bus stops &   Charging installment&BEB battery&& & & \\
         \midrule
                   0.2&5& 25 &   1.76&2.74&$4.504 \times 10^6$& $4.504 \times 10^6$& 0& 422\\
   &&  45 &   1.92&2.85&$4.774 \times 10^6$& $4.542\times 10^6$& 5.29\%&7200\\
           &25& 25 &   4.72&10.06&$14.788 \times 10^6$& $14.788 \times 10^6$& 0& 3595\\
   &&  45 &   4.24&10.65&$14.899\times 10^6$& $14.899\times 10^6$& 0&6224\\
           &45& 25 &   5.46&17.34&$22.808\times 10^6$& $22.808\times 10^6$& 0& 1048\\
   &&  45 &   4.96&18.63&$23.597\times 10^6$& $23.597\times 10^6$& 0&3695\\
   \cmidrule{1-9}
 0.8& 5& 25&   1.76&2.77&$4.539 \times 10^6$& $4.460\times 10^6$& 1.75\%&7200\\
 & & 45&   2.48&4.68&$7.166 \times 10^6$& $5.751\times 10^6$& 24.6\%&7200\\
           &25& 25&   4.64&10.68&$15.325\times 10^6$& $14.183 \times 10^6$& 3.8\%& 7200\\
   &&  45&   4&15.03&$19.030\times 10^6$& $14.183\times 10^6$& 25.5\%&7200\\
   &45&  25 &   5.92&17.72&$23.597\times 10^6$& $22.433\times 10^6$& 5.1\%&7200
\\
   &&  45 &   4.34&26.12&$30.469\times 10^6$& $21.115\times 10^6$& 44.3\%&7200\\
 \bottomrule
    \end{tabular}}
\end{table}

\textcolor{blue}{For high conservatism $(\theta=0.8)$, increasing network size leads to a decline in solution quality and a significant increase in computation times. Nevertheless, DRCC-CID remains capable of finding feasible solutions for large networks, though identifying the optimal solution takes longer (than the 2-hour time limit) due to the increased number of constraints and variables. However, these decisions are strategic and typically finalized before system operation begins, making computational times less critical.}

\textcolor{blue}{Electrification costs rise with network size and conservatism. The costs of BEB batteries account for over 60\% of the total electrification costs across all scenarios. To provide a tangible perspective on the charging infrastructure design and cost structure of CID, BoU-CID, and DRCC-CID models, the next section compares their solutions using three bus lines from the Rotterdam bus network.} 

\subsection{Deterministic and robust solutions on Rotterdam bus lines}
\label{sec:Rdam-case}
\textcolor{blue}{This section provides an in-depth examination of the CID, BoU-CID, and DRCC-CID solutions for the Rotterdam bus network including a detailed analysis of initial investment, charging infrastructure design and a review on long-term battery lifetime costs in \Cref{tab:breakdown_of_costs,tab:Rdam-results,tab:battery_life_time-results}. }

\textcolor{blue}{The total electrification costs of the CID, BoU-CID, and DRCC-CID models for different budget parameters are shown in \Cref{tab:breakdown_of_costs}. The first row of the table presents the deterministic CID model results, serving as a benchmark for the BoU-CID and DRCC-CID solutions listed below. Column 2 shows uncertainty levels, and Column 3 lists the total costs for each model. Column 4 shows the cost difference of solutions compared to the benchmark. Columns 5 and 6 present the percentage of the total electrification cost attributed to charging station installation and BEB battery capacity, respectively. The results show that modeling uncertainty using worst-case values (BoU-CID) leads to designs that are 67\% more expensive than the CID model under high conservatism, with increases observed in both charging station installation and battery capacity costs. In contrast, the DRCC-CID model incorporates additional data points into robust modeling, effectively managing cost increases while ensuring a robust design. Even under the more conservative setting, DRCC-CID incurs substantially lower costs compared to BoU-CID. The cost breakdown shows that DRCC-CID primarily addresses uncertainty by optimizing BEB battery capacity, reducing the need for costly charging station installations, and minimizing overall system costs.}
\begin{table}[h]
\centering
    \renewcommand{\arraystretch}{1.5} 
    \captionsetup{justification=raggedright,singlelinecheck=false,skip=2pt} 
    \centering

    \caption{\textcolor{blue}{Electrification costs breakdown of models}}
\tiny
\label{tab:breakdown_of_costs}

\textcolor{blue}{
    \begin{tabular}{p{1.5cm}p{2.5cm}p{1.5cm}p{1.5cm}p{2cm}p{1.5cm}}
    \toprule
        \textbf{Model}&\textbf{Uncertainty level}&\textbf{Total costs} (\euro)&\textbf{Relative difference} (\%)&\multicolumn{2}{l}{\textbf{Cost type percentage}}\\
        \cmidrule{5-6}& & & & CS installment& BEB battery\\
        \midrule
        \textbf{CID}& -& $8.413\times10^5$& \textbf{benchmark}& 41\%& 59\%\\
        \midrule
        \multirow{2}{*}{\textbf{BoU-CID}}& low: ($\Gamma = 0.2$)& $1.18\times10^6$& 24\%& 43\%& 57\%\\
         \cmidrule{2-6}& high: ($\Gamma = 0.8$)& $1.40\times10^6$& 67\%& 48\%& 52\%\\
         \midrule
                 \multirow{2}{*}{\textbf{DRCC-CID}}& low: ($\theta = 0.2$)& $9.81\times10^5$& 17\%& 27\%& 73\%\\
         \cmidrule{2-6}& high: ($\theta = 0.8$)& $9.98\times10^5$& 19\%& 26\%& 74\%\\
    \bottomrule
    \end{tabular}}
\end{table}

\textcolor{blue}{\Cref{tab:Rdam-results} highlights the charging infrastructure solutions for each model. Columns 3 and 4 display the charging station installation costs and the bus line number, respectively. Column 5 shows the Relative difference in battery capacity between the robust and deterministic solutions. The final two columns present the number and types of installed charging stations for each bus line. Detailed information on battery capacity values and charging station locations can be found in Appendix \cref{tab:BoU-performance,tab:DRCC-performance}.}

\textcolor{blue}{The deterministic CID solution installs five charging stations across the bus network. When accounting for worst-case energy consumption uncertainty in the BoU-CID solution, the required number of charging stations increases to seven or nine, depending on the allocated uncertainty budget. This reflects the trade-off between infrastructure expansion and uncertainty consideration. Increased charging station installations or the use of more powerful station types reduce battery capacity requirements, as charging becomes more frequent and efficient. For example, in the less conservative BoU-CID solution for bus line 40, battery capacity decreased compared to the CID solution due to the installation of additional charging stations, including fast-feeding types. However, allocating larger budgets to address uncertainty fluctuations results in simultaneous increases across all cost components, significantly driving up overall costs.}

\textcolor{blue}{DRCC-CID shows greater design stability compared to those relying solely on worst-case scenarios. In BoU-CID, charging station locations vary with the conservatism level, resulting in changing layouts and increased planning complexity. Conversely, DRCC-CID achieves a consistent and robust design for charging station locations and types, regardless of the conservatism degree. Additionally, DRCC-CID effectively manages larger uncertainty budgets by strategically increasing battery capacity in a controlled and cost-efficient manner.}
\begin{table}[h]
\centering
\renewcommand{\arraystretch}{1.5} 
\captionsetup{justification=raggedright,singlelinecheck=false,skip=2pt} 
\caption{\textcolor{blue}{Charging infrastructure solution of the models}}
\tiny
\label{tab:Rdam-results}
\textcolor{blue}{
\begin{tabular}{p{1.5cm}p{2cm}p{1cm}p{1cm}p{1.5cm}p{1cm}p{0.75cm}}
\toprule
\textbf{Model} & \textbf{Uncertainty level} & \textbf{Charging station cost} (\euro)& \textbf{Bus line number} & \textbf{Relative difference in battery capacity} & \multicolumn{2}{l}{\textbf{Number of Charging stations}} \\
\cmidrule{6-7}
& & & & & \textbf{Dedicated} & \textbf{Shared} \\
\midrule
\multirow{3}{*}{\textbf{CID}} 
& - & 340000 & 38 & \textbf{Benchmark} & 1FF & ~ \\
&   &        & 33 & ~                 & 1FF & 1FF \\
&   &        & 40 & ~                 & 1SS, 1FF & ~ \\
\midrule
\multirow{6}{*}{\textbf{BoU-CID}} 
& Low: ($\Gamma = 0.2$) & 440000 & 38 & +38\% & 1FF & ~ \\
&                       &        & 33 & +44\% & 1SS, 1FF & 1SS, 1FF \\
&                       &        & 40 & -7\%  & 2FF & ~ \\
\cmidrule{2-7}
& High: ($\Gamma = 0.8$) & 660000 & 38 & +35\% & 2FF & ~ \\
&                        &        & 33 & +79\% & 1SS, 2FF & 2FF \\
&                        &        & 40 & +42\% & 2FF & ~ \\
\midrule
\multirow{6}{*}{\textbf{DRCC-CID}} 
& Low: ($\theta = 0.2$) & 260000 & 38 & +80\% & 1FF & ~ \\
&                       &        & 33 & +20\% & 1SS & 1FF \\
&                       &        & 40 & +31\% & 1FF & ~ \\
\cmidrule{2-7}
& High: ($\theta = 0.8$) & 260000 & 38 & +85\% & 1FF & ~ \\
&                        &        & 33 & +27\% & 1SS & 1FF \\
&                        &        & 40 & +31\% & 1FF & ~ \\
\bottomrule
\end{tabular}}
Legend: \textbf{FF}: fast-feeding charging station; \textbf{SS}: Standard charging station
\end{table}

\textcolor{blue}{The charging infrastructure design and battery capacity significantly impact the BEB battery lifetime that is shown in \Cref{tab:battery_life_time-results}. The first three columns follow the same definitions as in the previous table. The fourth column indicates the optimal BEB battery capacity of each bus line. The fifth column presents the number of cycles to BEB battery failure, $N_{cycle}$, calculated using  \Cref{eq:n-cycle}. The cost per cycle is derived from \Cref{eq:cost-per-cyc}, and the sixth column highlights the difference in cost per cycle between the robust designs and the deterministic solution. Finally, the last column shows the average relative cost per cycle across the entire network. The results show that, although BoU-CID achieves lower costs per cycle for certain bus lines compared to the deterministic solution, it increases average costs across the entire network. In contrast, DRCC-CID consistently delivers lower costs per cycle for all bus lines, resulting in over 48\% cost savings per BEB battery cycle. This highlights the superior cost-efficiency of the DRCC-CID model in managing long-term battery lifetime expenses.}
\begin{table}[h]
\centering
    \renewcommand{\arraystretch}{1.5} 
    \captionsetup{justification=raggedright,singlelinecheck=false,skip=2pt} 
    \centering

    \caption{\textcolor{blue}{Battery lifetime of each model}}
\tiny
\label{tab:battery_life_time-results}
\textcolor{blue}{
    \begin{tabular}{p{1.5cm}p{2.5cm}p{1cm}p{1.5cm}p{1.5cm}p{1.5cm}p{2cm}p{3cm}}
        \toprule
        \textbf{Model}&\textbf{Uncertainty level}&\textbf{Bus line number}&\textbf{Battery capacity (kWh)} &\textbf{$N_{cycle}$}&{\textbf{Relative Cost per cycle}}& \textbf{Average cost difference for network}\\
        \midrule
        \multirow{3}{*}{\textbf{CID}}& -& 38 & 7.5 & 489940.95
&\textbf{benchmark}. &\\
         & & 33& 9.26 & 444205.1
& &\\
         & & 40 & 11.87 & 481227.57
& &\\
        \midrule
        \multirow{6}{*}{\textbf{BoU-CID}} & low: ($\Gamma = 0.2$)& 38 & 10.81 & 538176.8
& 19.8\%&\\
         & & 33& 12.75 & 620035.04
& 8.8\%&\\
         & & 40 & 11.09 & 525469.46
& -14.44\%&Avg: \textbf{+4.48}\%\\
         \cmidrule{2-7}& high: ($\Gamma = 0.8$)& 38 & 13.42 & 597594.36
& 62.74\%&\\
         & & 33& 12.49 & 741116.75
& -10.76\%&\\
         & & 40 & 16.82 & 704672.65
& -26.24\%&Avg.: \textbf{+8.58}\%\\
         \midrule
                 \multirow{6}{*}{\textbf{DRCC-CID}} & low: ($\theta = 0.2$)& 38 & 9.01 & 1104716.88
& -51.7\%&\\
         & & 33& 16.63 & 1506432.1
& -41.6\%&\\
         & & 40 & 15.54 & 1406003.3
& -55.19\%&Avg.:\textbf{ -49.5}\%\\
         \cmidrule{2-7}& high: ($\theta = 0.8$)& 38 & 9.51 & 1104716.88
& -49.01\%&\\
         & & 33& 17.13 & 1506432.1
& -39.84\%&\\
         & & 40 & 15.54 & 1406003.3& -55.19\%&Avg.:\textbf{ -48}\%\\
        \bottomrule
    \end{tabular}}
\end{table}

\subsection{Monte Carlo simulation: Feasibility analysis of the results of the three models}
\label{sec:simulation-monte}
\textcolor{blue}{Comparing the optimal solution costs and performance of robust and deterministic models can be misleading since each model is optimized for a different dataset. Specifically, the deterministic CID model relies on nominal realizations of energy consumption data, the BoU-CID model addresses worst-case scenarios, and the DRCC-CID model integrates sparse data into robust modeling. To ensure a fair comparison, we adopt the approach of \citep{avishan2023electric}, employing a Monte Carlo simulation experiment to evaluate the optimality and feasibility of the models on equal terms. In this method, energy consumption data are randomly sampled from various distributions. A charging infrastructure design is deemed infeasible for a given scenario if the battery level upon arrival at a stop falls below the minimum allowed limit, indicating that the design cannot support the trip under the realized conditions. This process is repeated for each design solution, and the number of infeasible scenarios is recorded based on the sampled energy data. The entire procedure is repeated 100 times, and the feasibility rates for each bus line and the average feasibility rate across the network are reported. The pseudo-code for the simulation algorithm is presented in \ref{simulation_pseudo}.}

\textcolor{blue}{\textbf{In-sampling performance of the models.} We evaluate the effectiveness of the robust approach in addressing uncertainty by employing the Monte Carlo simulation method described earlier. Uncertainty parameters are sampled from four different distributions: uniform ($\mathcal{U}$) and triangular ($\mathcal{T}$), as well as asymmetric right-triangular and left-triangular distributions. The robust models are tested under 20\% and 80\% conservatism levels, with results summarized in \Cref{tab:in-sample-feas-results}. The first and second columns specify the simulation number and the sampling distribution for the uncertainty parameters, respectively. The third and fourth columns detail the model name and budget parameter settings. Column 5 reports the total electrification costs associated with each sampling distribution. Columns 6 to 8 present the feasibility percentages of the obtained designs over 100 simulation runs, and the final column provides the average feasibility performance across the entire bus network.}

\textcolor{blue}{The four sampling distributions are designed to evaluate model performance under varying data characteristics. The first distribution is a uniform distribution centered on the mean energy consumption, with a range extending to the maximum deviation. The second is a symmetric triangular distribution. The third and fourth are asymmetric right- and left-triangular distributions, representing scenarios where the maximum deviation from the mean has the least and highest probability of occurrence, respectively. This approach allows us to assess the impact of incorporating sample data into robust modeling, particularly when the range or average values are poor representations of the actual data.}
\begin{table}[h]
\centering
    \renewcommand{\arraystretch}{1.5} 
    \captionsetup{justification=raggedright,singlelinecheck=false,skip=2pt} 
    \centering

    \caption{\textcolor{blue}{In-sample simulation results for CID, BoU-CID, and DRCC-CID}}
\tiny
\label{tab:in-sample-feas-results}
\textcolor{blue}{
    \begin{tabular}{lp{1.75cm}p{1.5cm}p{2cm}p{1.5cm}p{0.5cm}p{0.5cm}p{0.5cm}p{1.5cm}}
        \toprule
          nr.&\textbf{Sampling distribution}&\multicolumn{2}{l}{\textbf{Uncertainty level of solved problem}}& \textbf{Total cost} (\euro)&\multicolumn{3}{l}{\textbf{Bus line number}}&\textbf{Average feasibility of the bus network}\\
          \cmidrule{6-8}
          && model& &  & 38& 33&40& \\
        \midrule
          1&$\mathcal{U}[0, \hat{\mu}]$&\textbf{CID}& -&  $8.1\times10^5$&56\%& 97\%&67\%& 73\%\\
         \cmidrule{3-9}&&\textbf{BoU-CID}& low: ($\Gamma = 0.2$)&  $8.76\times10^5$&100\%& 60\%&100\%& 87\%\\
           &&& high: ($\Gamma = 0.8$)&  $1.31\times10^6$&\multicolumn{4}{l}{Fully feasible}\\
           \cmidrule{3-9}&&\textbf{DRCC-CID}& low: ($\theta = 0.2$)&  $8.1\times10^5$&51\%& 93\%&92\%& 79\%\\
           &&& high: ($\theta = 0.8$)&  $8.31\times10^5$&96\%& 93\%&97\%& 95\%\\
        \midrule
         2&$\mathcal{T}[0, \hat{\mu}, \frac{\hat{\mu}}{2}]$&\textbf{CID}& -&  $7.56\times10^5$&52\%& 68\%&70\%& 63\%\\
         \cmidrule{3-9}&&\textbf{BoU-CID}& low: ($\Gamma = 0.2$)&  $9.58\times10^5$&\multicolumn{4}{l}{Fully feasible}\\
           &&& high: ($\Gamma = 0.8$)&  $1.32\times10^6$&\multicolumn{4}{l}{Fully feasible}\\
        \cmidrule{3-9}&&\textbf{DRCC-CID}& low: ($\theta = 0.2$)&  $8.72\times10^5$&43\%& 92\%&92\%& 76\%\\
           &&& high: ($\theta = 0.8$)&  $8.75\times10^5$&96\%& 93\%&97\%& 93\%\\
        \midrule
          3&$\mathcal{T}[0, \hat{\mu},0]$&\textbf{CID}& -&  $6.32\times10^5$&68\%& 71\%&76\%& 72\%\\
         \cmidrule{3-9}&&\textbf{BoU-CID}& low: ($\Gamma = 0.2$)&  $7.51\times10^5$&100\%& 98\%&100\%& 99\%\\
           &&& high: ($\Gamma = 0.8$)&  $1.24\times10^6$&\multicolumn{4}{l}{Fully feasible}\\
        \cmidrule{3-9}&&\textbf{DRCC-CID}& low: ($\theta = 0.2$)&  $6.43\times10^5$&74\%& 99\%&91\%& 88\%\\
           &&& high: ($\theta = 0.8$)&  $6.63\times10^5$&99\%& 99\%&95\%& 98\%\\
        \midrule
         4&$\mathcal{T}[0, \hat{\mu}, \hat{\mu}]$&\textbf{CID}& -&  $8.94\times10^5$&40\%& 20\%&75\%& 45\%\\
         \cmidrule{3-9}&&\textbf{BoU-CID}& low: ($\Gamma = 0.2$)&  $1.29\times10^6$&95\%& 98\%&98\%& 97\%\\
           &&& high: ($\Gamma = 0.8$)&  $1.38\times10^6$&\multicolumn{4}{l}{Fully feasible}\\
        \cmidrule{3-9}&&\textbf{DRCC-CID}& low: ($\theta = 0.2$)&  $9.19\times10^5$&36\%& 80\%&91\%& 69\%\\
           &&& high: ($\theta = 0.8$)&  $9.97\times10^5$&65\%& 92\%&89\%& 82\%\\
        \midrule
    \end{tabular}}
\end{table}

\textcolor{blue}{The results demonstrate that while the deterministic CID model performs well when energy consumption data is centered around the average (instances 1 and 2) or skewed towards lower values (instance 3), it fails to maintain feasibility in more extreme scenarios, such as instance 4, where higher-than-average values dominate, leading to a 55\% infeasibility rate. As expected, BoU-CID ensures full feasibility across almost all scenarios, particularly under higher uncertainty budgets. However, the DRCC-CID model strikes a balance by significantly improving feasibility over CID, especially when average values poorly represent the data distribution, while consistently achieving lower costs than BoU-CID. This highlights DRCC-CID's ability to provide robust and cost-efficient solutions under diverse uncertainty conditions.}

\textcolor{blue}{\textbf{Out-of-sampling performance of the models.} This section explores the performance of the models under unobserved conditions by altering the probability distributions used for sampling. In this modified simulation experiment, the sampling distributions for testing differ from the distribution utilized in designing the models, differing in their supports and distributional characteristics. \Cref{tab:out_of-sample-feas-results} summarizes the findings, with the second column indicating the distribution used to achieve the design decisions of the model and the third column is the sampling distribution used for testing the performance of the model in previously unconsidered scenarios. The fifth column shows the relative difference in total electrification costs for robust models compared to the CID solution.} 

\textcolor{blue}{The first and second simulation numbers examine cases where the test energy consumption values can exceed the design data by 20\% in the range. Simulations 3 and 4 explore scenarios where the test and design distributions differ in their characteristics: a symmetric triangular design distribution is compared against a left-triangular distribution and a symmetric triangular distribution with a shifted mode. These test distributions assign the highest and lowest probabilities, respectively, to maximum deviations from the average energy consumption.}
\begin{table}[h]
\centering
    \renewcommand{\arraystretch}{1.5} 
    \captionsetup{justification=raggedright,singlelinecheck=false,skip=2pt} 
    \centering

    \caption{\textcolor{blue}{Out-of-sample simulation results for CID, BoU-CID, and DRCC-CID}}
\tiny
\label{tab:out_of-sample-feas-results}
\textcolor{blue}{
    \begin{tabular}{lp{1.25cm}p{1.6cm}p{1.5cm}p{1.75cm}p{1cm}p{0.5cm}p{0.5cm}p{0.5cm}p{1.5cm}}
        \toprule
           nr.&\textbf{Design distribution}&\textbf{Out-of-sampling distribution}&\multicolumn{2}{l}{\textbf{Uncertainty level of solved problem}} &\textbf{Total cost}&\multicolumn{3}{l}{\textbf{Bus line number}}&\textbf{Average feasibility of the bus network}\\
           \cmidrule{6-9}&&& model&  &Relative difference& 38& 33&40& \\
        \midrule
           1&$\mathcal{U}[0, \hat{\mu}]$&$\mathcal{U}[0, \textbf{1.2}\hat{\mu}]$&\textbf{CID}& - &bench.&12\%& 17\%&14\%& 14\%\\
           \cmidrule{4-10}&&&\textbf{BoU-CID}& low: ($\Gamma = 0.2$) &+8\%&68\%& 59\%&40\%& 56\%\\
            &&&& high: ($\Gamma = 0.8$) &+62\%&\multicolumn{4}{l}{Fully feasible}\\
            \cmidrule{4-10}&&&\textbf{DRCC-CID}& low: ($\theta = 0.2$) &+0.2\%&14\%& 48\%&49\%& 37\%\\
            &&&& high: ($\theta = 0.8$) &+3\%&53\%& 50\%&65\%& 56\%\\
        \midrule
          2&$\mathcal{T}[0, \hat{\mu}, \hat{\mu}]$&$\mathcal{T}[0, \textbf{1.2}\hat{\mu}, 
\textbf{1.2} \hat{\mu}]$&\textbf{CID}& - &bench.&2\%& 4\%&4\%& 3\%\\
          \cmidrule{4-10}&&&\textbf{BoU-CID}& low: ($\Gamma = 0.2$) &+44\%&49\%& 40\%&35\%& 41\%\\
            &&&& high: ($\Gamma = 0.8$) &+54\%&\multicolumn{4}{l}{Fully feasible}\\
         \cmidrule{4-10}&&&\textbf{DRCC-CID}& low: ($\theta = 0.2$) &+3\%&4\%& 20\%&20\%& 15\%\\
            &&&& high: ($\theta = 0.8$) &+12\%&39\%& 54\%&48\%& 47\%\\
        \midrule
           3&$\mathcal{T}[0, \hat{\mu}, \frac{\hat{\mu}}{2}]$&$\mathcal{T}[\bar{\mu}, \hat{\mu}, \hat{\mu}]$&\textbf{CID}& - &bench.&1\%& 0\%&0\%& 0\%\\
          \cmidrule{4-10}&&&\textbf{BoU-CID}& low: ($\Gamma = 0.2$) &+27\%&49\%& 58\%&32\%& 46\%\\
            &&&& high: ($\Gamma = 0.8$) &+75\%&\multicolumn{4}{l}{Fully feasible}\\
         \cmidrule{4-10}&&&\textbf{DRCC-CID}& low: ($\theta = 0.2$) &+15\%&4\%& 20\%&20\%& 15\%\\
            &&&& high: ($\theta = 0.8$) &+16\%&17\%& 20\%&13\%& 17\%\\
            \midrule
 4& $\mathcal{T}[0, \hat{\mu}, \frac{\hat{\mu}}{2}]$& $\mathcal{T}[0, \hat{\mu}, \textbf{0.8}\bar{\mu}]$& \textbf{CID}& - &bench.& 14\%& 18\%& 19\%&17\%\\
 \cmidrule{4-10}& & & \textbf{BoU-CID}& low: ($\Gamma = 0.2$) &+27\%& 96\%& 89\%& 82\%&89\%\\
 & & & & high: ($\Gamma = 0.8$) &+75\%& \multicolumn{4}{l}{Fully feasible}\\
 \cmidrule{4-10}& & & \textbf{DRCC-CID}& low: ($\theta = 0.2$) &+15\%& 37\%& 41\%& 43\%&40\%\\
 & & & & high: ($\theta = 0.8$) &+16\%& 63\%& 71\%& 81\%&71\%\\
 \bottomrule
    \end{tabular}} 
\end{table}

\textcolor{blue}{The results reveal that relying solely on average energy consumption for the design of the charging infrastructure, as in the CID model, leads to high infeasibility in all scenarios, with average feasibility rates ranging from 0\% to 17\%. The CID model is particularly inadequate in the third simulation, where the probability of worst-case values is high, leading to complete infeasibility. By accounting for uncertainty, BoU-CID ensures full feasibility under high conservatism (80\%). However, under low conservatism (20\%), feasibility drops to 41\% in the second simulation, where both the range and data characteristics deviate from the design.The DRCC-CID model consistently outperforms the deterministic CID model in all simulations, demonstrating improved feasibility under uncertainty. When test data differ only in range but retain similar distributional properties (simulation numbers 1 and 2), DRCC-CID maintains a reasonable level of feasibility while avoiding the excessive costs associated with BoU-CID. However, in Simulation 3, where test data are centered around worst-case values, DRCC-CID's feasibility decreases. Although BoU-CID ensures full feasibility in such extreme scenarios, it does so at significantly higher costs due to its focus on the worst case.}

\textcolor{blue}{The DRCC-CID model effectively integrates empirical data into the robust optimization framework, excelling when worst-case scenarios are rare, and data are skewed away from extremes, as demonstrated in Simulation number 4. DRCC-CID minimizes unnecessary costs and avoids overcompensation for unlikely extremes, striking an optimal balance between feasibility and cost-efficiency. This data-driven approach ensures robust performance and resource efficiency in practical applications, making it a compelling and cost-efficient alternative to worst-case-focused models like BoU-CID.}

\subsubsection{Effect of data quantity on data-driven distributionally robust optimization solutions}
\label{sec:quantity-drcc}
\begin{figure}[ht]
    \centering
    \begin{subfigure}[b]{0.48\textwidth}
        \centering
        \includegraphics[width=\textwidth]{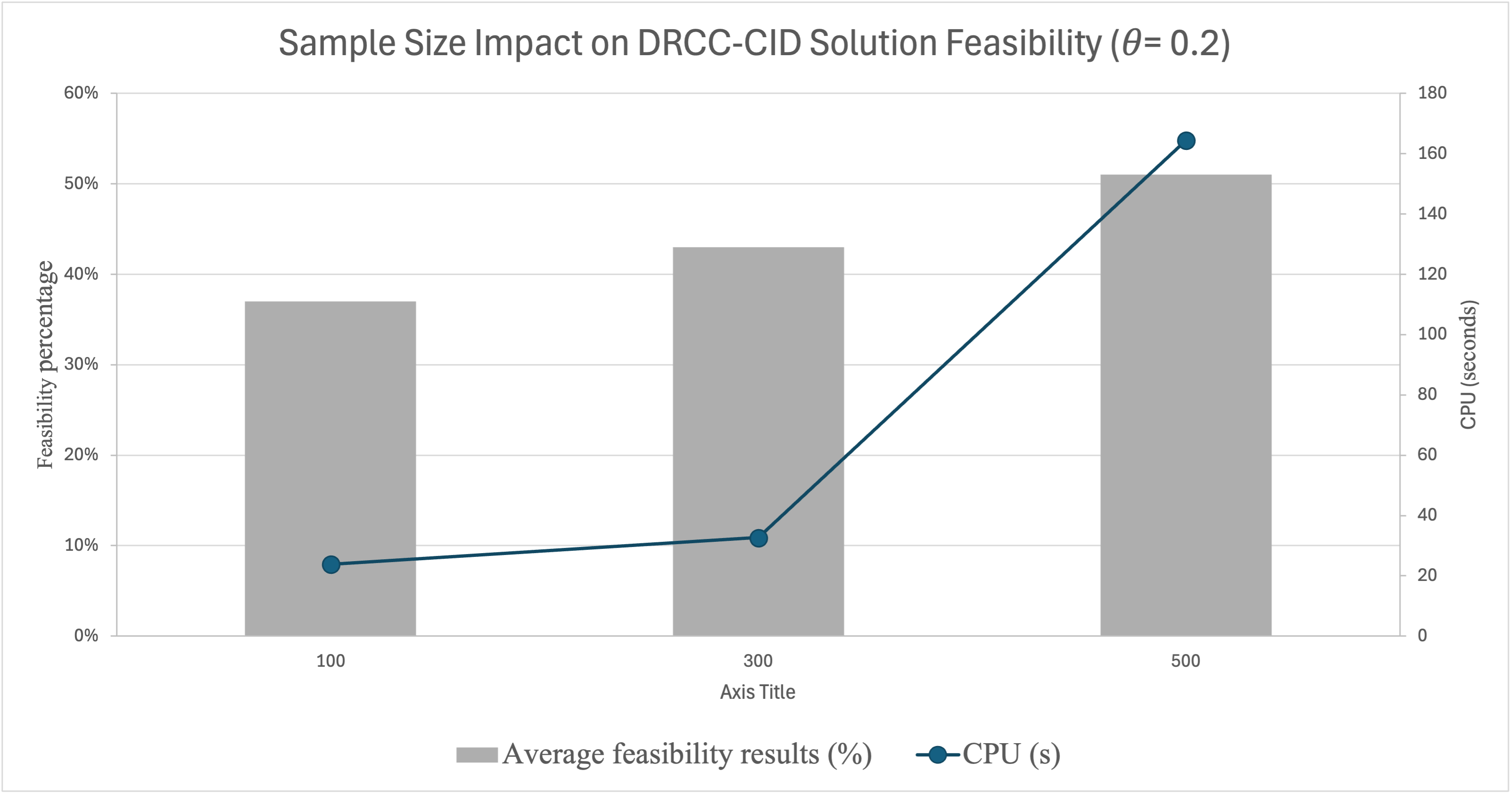}
        \caption{DRCC-CID with 20\% budget parameter}
        \label{fig:DataQ-0.2}
    \end{subfigure}
    \hfill
    \begin{subfigure}[b]{0.48\textwidth}
        \centering
        \includegraphics[width=\textwidth]{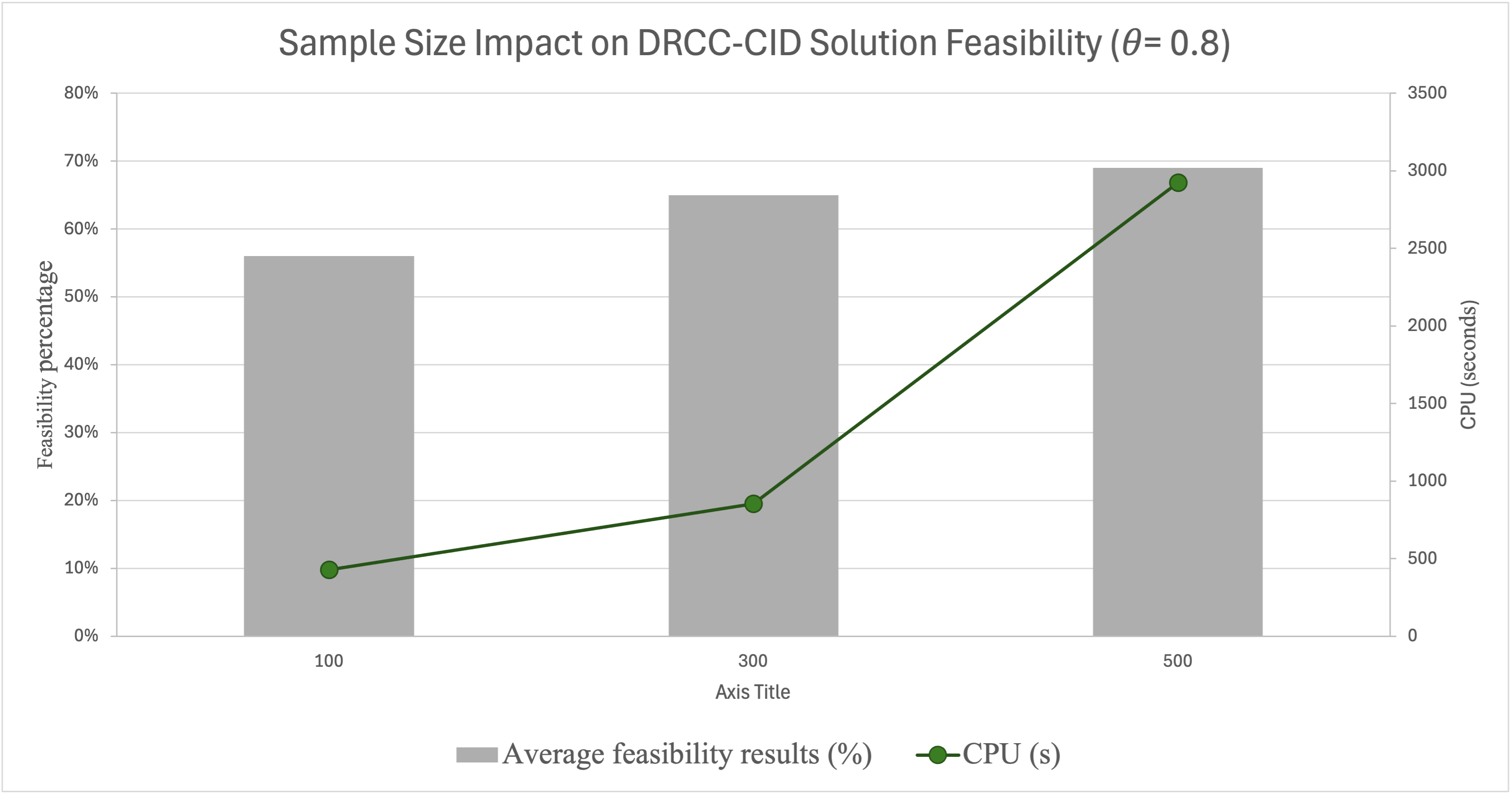}
        \caption{DRCC-CID with 80\% budget parameter}
        \label{fig:dataQ-0.8}
    \end{subfigure}
    \caption{\textcolor{blue}{Impact of sample size on feasibility and computational time of the DRCC-CID model}}
    \label{fig:DataQ-feasibility-CPU}
\end{figure}

\textcolor{blue}{The DRCC-CID model is designed to perform effectively with sparse sample data. In this section, we evaluate the impact of sample size on the DRCC-CID performance, focusing on the feasibility of the design and computational time. In the previous section, \Cref{tab:out_of-sample-feas-results}, the out-of-sample performance was evaluated using a baseline of 100 energy consumption data points between consecutive stops. Here, we extend the analysis by solving the DRCC-CID model with varying sample sizes of 100, 300, and 500 data points. This analysis uses the first simulation scenario from \Cref{tab:out_of-sample-feas-results} as the design and test distributions. Detailed results for this and the third simulation scenario are provided in \ref{appendix:quantity}. The findings for the first simulation scenario are summarized in \Cref{fig:DataQ-feasibility-CPU}, illustrating the impact of sample size on design feasibility and computational efficiency. Results are presented for two budget parameters: 20\% (\ref{fig:DataQ-0.2}) and 80\% (\ref{fig:dataQ-0.8}).
}

\textcolor{blue}{For a 20\% budget parameter, \Cref{fig:DataQ-0.2}, increasing the sample size from 100 to 500 data points improves feasibility from 37\% to 51\%, while CPU time increases significantly from 24 seconds to nearly 15 minutes. Similarly, for an 80\% budget parameter, \Cref{fig:dataQ-0.8}, the feasibility improves from 56\% to 69\% as the sample size increases from 100 to 500 data points, but the computational time rises steeply from 7 minutes to 48 minutes. These results indicate that while increasing sample size enhances feasibility, the improvement is relatively modest compared to the computational time.}

\textcolor{blue}{A sensitivity analysis of the budget parameters reveals that the feasibility of the design is more strongly influenced by the parameters of the robust model, i.e. budget parameter, than by the sample size. This underscores the effectiveness of the DRCC-CID model in maintaining performance even with sparse data, making it a practical and computationally efficient choice in settings with limited data or resource constraints.}

\subsection{Sensitivity analysis in charging infrastructure design}
\label{sec:sensitivity}
\textcolor{blue}{In this section, we examine the sensitivity of the models to key economic factors, including charging station installation costs and BEB fleet size. Although the BEB fleet size is treated as a known parameter in our case study, we investigate how variations in fleet size influence the cost structure of the models, offering valuable insights into their adaptability and robustness under different economic scenarios.}

\subsubsection{Economic factor sensitivity}
\textcolor{blue}{This section investigates the sensitivity of the models to charging station installation costs and their influence on total onboard battery capacities. The analysis is conducted by fixing the budget parameters for both robust modeling approaches at 20\% and 80\%, representing low and high conservatism levels, respectively.}

\begin{figure}[htbp]
    \centering
    \begin{subfigure}[b]{0.45\textwidth}
        \centering
        \includegraphics[width=\textwidth]{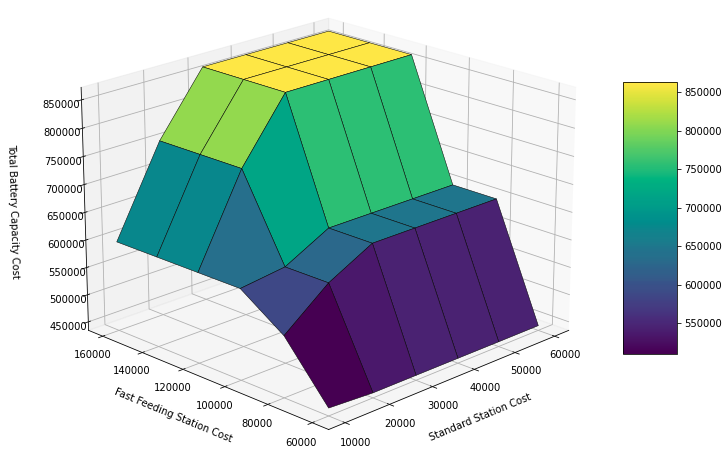}
        \caption{BoU-CID solution with $\Gamma$=0.2}
        \label{fig:BoU_0.2_cost_sens}
    \end{subfigure}
    \hfill
     \begin{subfigure}[b]{0.45\textwidth}
        \centering
        \includegraphics[width=\textwidth]{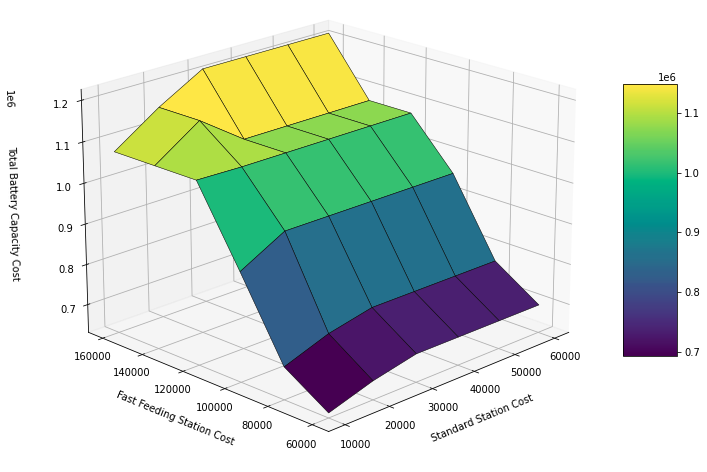}
        \caption{BoU-CID solution with $\Gamma$=0.8}
        \label{fig:BoU-0.8_cost_sens}
    \end{subfigure}
     \hfill
    \begin{subfigure}[b]{0.45\textwidth}
        \centering
        \includegraphics[width=\textwidth]{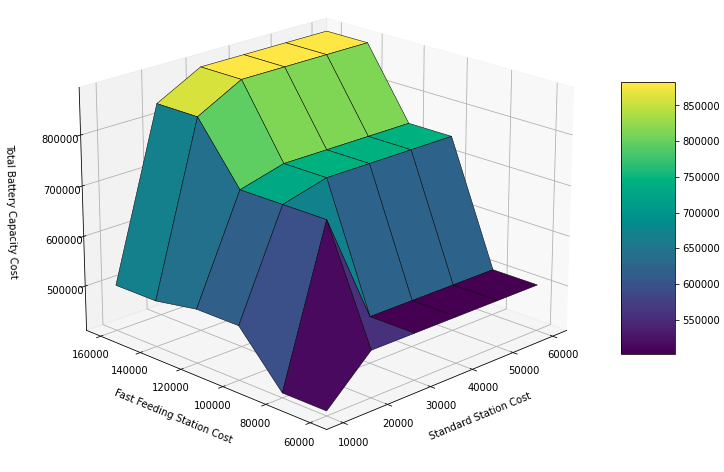}
        \caption{DRCC-CID solution with $\theta = 0.2$}
        \label{fig:drcc-cost-sens20_sens}
    \end{subfigure}
    \hfill
    \begin{subfigure}[b]{0.45\textwidth}
        \centering
        \includegraphics[width=\textwidth]{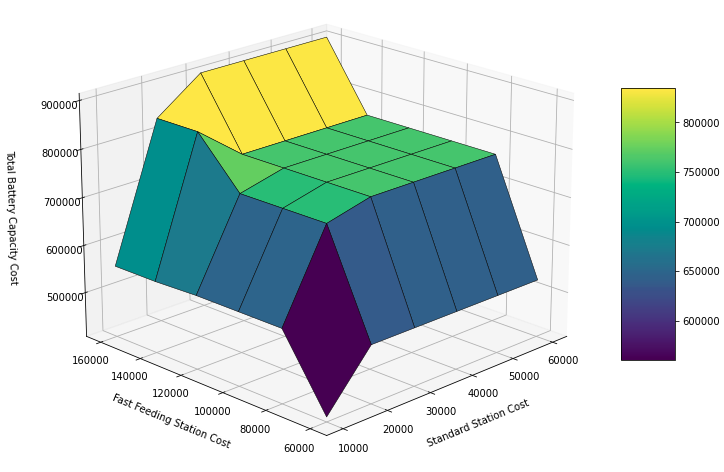}
        \caption{DRCC-CID solution with $\theta = 0.8$}
        \label{fig:drcc_cost-80_sens}
    \end{subfigure}
    \caption{\textcolor{blue}{The impact of charging station costs on total battery capacity of robust models}}
    \label{fig:drcc_sensitivity}
\end{figure}

\textcolor{blue}{For the analysis, we assume that the installation cost of standard chargers ranges from $\alpha_{SS}=[10000,60000]$ euros, and for fast-feeding chargers, $\alpha_{FF}=[60000,160000]$} euros. The maximum cost values are based on \cite{lajunen2014energy}. \textcolor{blue}{\Cref{fig:drcc_sensitivity} illustrates the impact of variations in charging station costs on total battery capacity of BoU-CID and DRCC-CID models, in the first and second rows respectively.} 

\textcolor{blue}{At a low conservatism level, the BoU-CID  \Cref{fig:BoU_0.2_cost_sens} increases total battery capacity more significantly than DRCC-CID  \ref{fig:drcc-cost-sens20_sens}. However, the BoU-CID stabilizes total battery capacity at the highest charging station costs, indicating a preference for installing charging stations rather than addressing uncertainty through further capacity increases. In contrast, the DRCC-CID continues to rely more heavily on capacity adjustments, showing less stability at higher charging station costs.}

\textcolor{blue}{At high conservatism levels, both robust models prioritize installing more fast-feeding charging stations and show indifference to increases in standard charging station costs. The BoU-CID solution in \Cref{fig:BoU-0.8_cost_sens} consistently increases total battery capacity while installing more stations.  In comparison, the DRCC-CID solution \Cref{fig:drcc_cost-80_sens}, demonstrates greater stability in total battery capacity, emphasizing cost-efficient solutions rather than excessive capacity increases. }

\subsubsection{Fleet size sensitivity}
\textcolor{blue}{This section evaluates the impact of BEB fleet size on charging infrastructure design and the cost structure of the solution. Testing different fleet sizes is essential to understanding how operational decisions, such as increasing fleet size for improved reliability or contingency planning, affect the overall electrification strategy. While the base scenario, assuming 10 BEBs per bus line, aligns with the case study assumptions, scenarios with 5 and 15 BEBs per line allow us to explore the trade-offs between fleet size and other factors, such as charging station installations and battery capacity.}

\textcolor{blue}{The results, summarized in \Cref{tab:fleet-size-sens}, are organized as follows: the first and second columns specify the model name and the corresponding uncertainty level, respectively. Column 3 presents the fleet sizes analyzed (5, 10, and 15 BEBs per bus line), while Column 4 details the relative differences in charging station installation costs and battery capacity compared to the benchmark. Finally, the last columns provide the cost composition, showing the percentage contributions of charging station installations and battery capacity to the total electrification costs for each solution.}

\begin{table}[h]
\centering
    \renewcommand{\arraystretch}{1.5} 
    \captionsetup{justification=raggedright,singlelinecheck=false,skip=2pt} 
    \centering

    \caption{\textcolor{blue}{Impact of BEB fleet size on charging station installation and battery capacity costs}}
\tiny
\label{tab:fleet-size-sens}

\textcolor{blue}{
    \begin{tabular}{p{1.5cm}p{2.5cm}p{1cm}p{1.75cm}p{1.25cm}p{1cm}p{1cm}}
    \toprule
        \textbf{Model}&\textbf{Uncertainty level}& \textbf{BEB fleet size}&\multicolumn{2}{l}{\textbf{Relative difference}}&\multicolumn{2}{l}{\textbf{Cost type percentage}}\\
        \cmidrule{4-7}& &  &CS installation costs (\euro)& Battery capacity (kWh)& CS installment (\euro)& Battery capacity (\euro)\\
        \midrule
        \textbf{CID}& -&  10&\textbf{benchmark}& & & \\
 & & 5& -& +18\%& 0\%&100\%\\
 & & 15& +17\%& -58\%& 64\%&36\%\\
        \midrule
        \multirow{2}{*}{\textbf{BoU-CID}}& low: ($\Gamma = 0.2$)&  10&\textbf{benchmark}& & & \\
 & & 5& -69\%& +67\%& 29\%&71\%\\
 & & 15& +0.09\%& +177\%& 42\%&58\%\\
         \cmidrule{2-7}& high: ($\Gamma = 0.8$)&  10&\textbf{benchmark}& & & \\
 & & 5& -50\%& +58\%& 37\%&63\%\\
 & & 15& +18\%& -56\%& 47\%&53\%\\
         \midrule
                 \multirow{2}{*}{\textbf{DRCC-CID}}& low: ($\theta = 0.2$)&  10&\textbf{benchmark}& & & \\
 & & 5& -39\%& +22\%& 27\%&73\%\\
 & & 15& +5\%& -30\%&39\% &61\%\\
         \cmidrule{2-7}& high: ($\theta = 0.8$)&  10&\textbf{benchmark}& & & \\
 & & 5& -30\%& +24\%& 25\%&75\%\\
 & & 15& +7\%& -26\%& 39\% &61\%\\
 \bottomrule
    \end{tabular}}
\end{table}

\textcolor{blue}{In the deterministic CID model, a reduction in fleet size to 5 results in no charging station installations and an 18\% increase in battery capacity compared to the base solution. Conversely, for a fleet size of 15, the model opts to install charging stations and reduce battery capacity, effectively balancing total electrification costs. A similar pattern of adjusting charging station installations and battery capacity with fleet size is observed in the robust models.}

\textcolor{blue}{However, an exception occurs in the BoU-CID model under a 20\% conservatism level. For a fleet size of 15, the model retains a similar charging station installation design as for a fleet size of 10, while increasing battery capacity significantly to hedge against uncertainty. In contrast, DRCC-CID demonstrates a more balanced approach across both 20\% and 80\% budget parameters. For a fleet size of 5, DRCC-CID reduces charging station installation costs by 30-39\% and increases battery capacity by only 22-24\%, minimizing total costs. For a fleet size of 15, it increases charging station costs moderately by 5-7\% and reduces battery capacity by 26-30\%, achieving a more cost-efficient and balanced design compared to BoU-CID.}

\textcolor{blue}{A consistent trend emerges across all models. With larger BEB fleet sizes, the models focus on installing additional charging stations to accommodate increased demand, thereby reducing reliance on onboard battery capacity. For smaller fleet sizes, charging station installations are minimized, and the models compensate by increasing battery capacity to ensure feasibility. It is important to note that fleet size decisions are closely linked to bus scheduling considerations, which are beyond the scope of this study.}

\section{Conclusion}
\label{sec:conclusion}
\textcolor{blue}{In this study, we investigate the value of energy consumption data collection and its impact on the cost efficiency of robust Battery Electric Bus (BEB) charging infrastructure design. To address the challenges posed by sparse and uncertain energy consumption data, we implement two robust modeling approaches. The first approach, range-Bound scenario planning, uses a box uncertainty set with a budget parameter to manage conservatism, relying on the range characteristics of the random parameter. The second approach, a data-driven distributionally robust optimization, directly incorporates sparse sample data into the robust optimization process. Our analysis evaluates the effectiveness of these uncertainty modeling techniques, with a focus on cost efficiency and solution feasibility.}

\textcolor{blue}{We validate these approaches using hypothetical grid networks and a real-world case study of the Rotterdam bus network. Both robust models are capable of delivering optimal solutions, even for large-scale bus networks. However, the results show that relying solely on the range of uncertain energy consumption values can lead to significantly inflated electrification costs. Conversely, the data-driven distributionally robust optimization consistently demonstrates superior cost-efficiency and high feasibility, particularly when worst-case scenarios are infrequent, and data are skewed away from extreme values. Sensitivity analyses on sample size and computation time reveal that the data-driven distributionally robust optimization model maintains robust performance even with limited data, making it a practical and computationally efficient option in resource-constrained settings.}

\textcolor{blue}{This study highlights the substantial influence of modeling approaches on BEB infrastructure design and electrification costs. Future research should explore additional dimensions of uncertainty, such as supply-side factors, and extend to tactical and operational decision-making, including the integration of bus schedules with BEB fleet sizing. While this study assumes a complete transition to BEBs, future work could examine the implications of partial electrification and the associated infrastructure adjustments. Addressing these aspects will further enhance the robustness and efficiency of BEB infrastructure design, contributing to more sustainable and reliable urban transit systems.} 
 
\bibliographystyle{elsarticle-harv}
\bibliography{J1}

\newpage
\appendix

\section{Hypothetical grid network}
\begin{figure}[h!]
    \centering
    \begin{subfigure}[b]{0.3\textwidth}
        \centering
        \includegraphics[width=\textwidth]{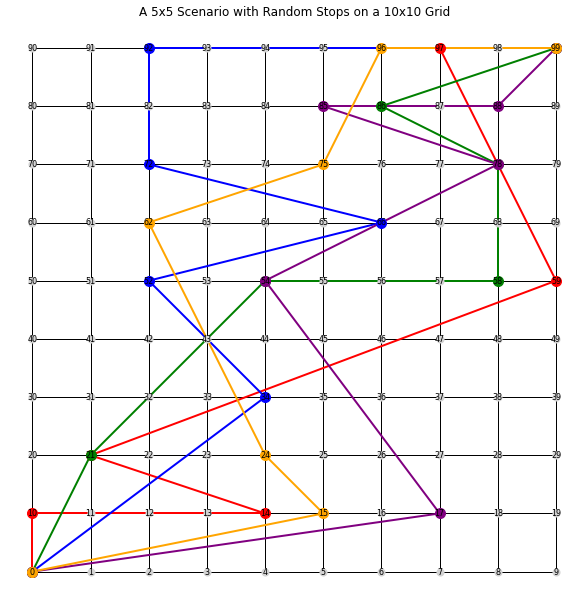}
        \caption{Scenario of 5 bus lines, each consisting of 5 bus stops}
        \label{fig:5-5grid}
    \end{subfigure}
    \hfill
    \begin{subfigure}[b]{0.3\textwidth}
        \centering
        \includegraphics[width=\textwidth]{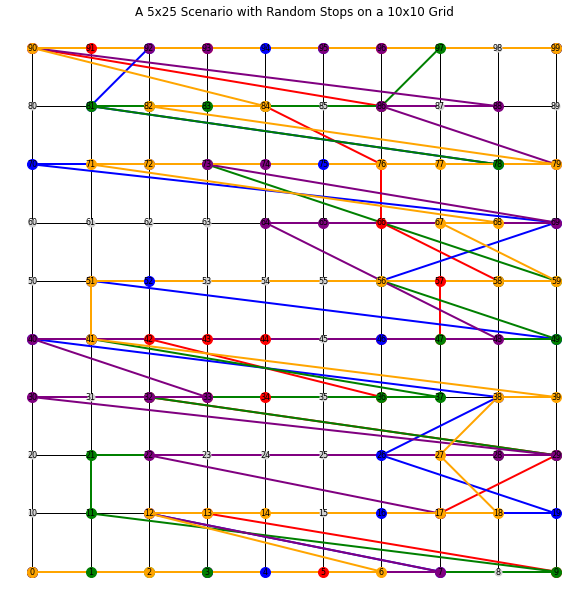}
        \caption{5 bus lines, each containing 25 bus stops}
        \label{fig:5-25grid}
    \end{subfigure}
    \hfill
    \begin{subfigure}[b]{0.3\textwidth}
        \centering
        \includegraphics[width=\textwidth]{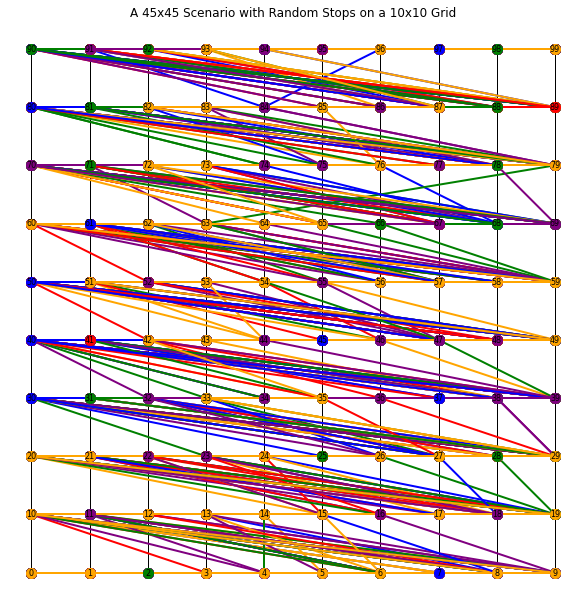}
        \caption{45 bus lines, each containing 45 bus stops}
        \label{fig:45-45grid}
    \end{subfigure}
    \caption{\textcolor{blue}{Diverse network configuration scenarios in hypothetic grid network with varying level of inter-connectivity}}
    \label{fig:hypo-grids}
\end{figure}

\textcolor{blue}{\Cref{fig:hypo-grids} presents a schematic view of generated bus networks over a $10\times 10$ grid. Each node on the grid represents a potential bus stop or charging station location. The scenarios depict varying densities and distances between stops, reflecting diverse network configurations. \Cref{fig:5-5grid} illustrates a sparse network with 5 bus lines, each containing 5 stops, characterized by widely spaced stops. \Cref{fig:5-25grid} represents a more dense network with more closely spaced stops. Finally, \Cref{fig:45-45grid} shows the largest network, consisting of 45 bus lines, each with 45 stops, representing a highly dense and interconnected bus network.}

\newpage
\section{Performance of BoU-CID and DRCC-CID models under varying conservatism levels}
\begin{table}[htbp]
    \centering
    \renewcommand{\arraystretch}{1.5} 
    \captionsetup{justification=raggedright,singlelinecheck=false,skip=2pt} 
    \caption{\textcolor{blue}{BoU-CID solution sensitivity in terms of budget parameter (\textbf{$\Gamma$})}}
    \tiny
    \textcolor{blue}{
    \begin{tabular}{lp{1.5cm}p{1.125cm}p{1.5cm}p{1cm}p{1cm}p{1cm}p{0.5cm}p{1.5cm}p{2.5cm}}
    \toprule
       \multirow{2}{*}{(\textbf{$\Gamma$})} &  \multicolumn{3}{l}{\textbf{Electrification cost}}&\multicolumn{3}{l}{\textbf{Battery capacity (kWh)}} & \multirow{2}{*}{\textbf{CPU}(s)}&\multicolumn{2}{l}{\textbf{Network configuration}}\\
       \cmidrule{2-7}
       \cmidrule{9-10}
       &OFV (\euro)& CS installation&Relative difference&Line 33 &Line 38 &Line 40& & SS&FF\\ 
       \midrule
        0 & $8.413\times10^5$& 340000&\textbf{benchmark}.&9.26&7.5&11.87&9.41&\parbox[t]{1.5cm}{\{Delft, Balthasar Van Der Polweg\} }&\parbox[t]{3cm}{\{Paradijsplein, Schipluiden- Halfwege, return-De Lugt, 2e Hogenbanweg\}}\\
        0.2 &$1.18\times10^6$& 440000
&+24\%&12.75&10.81&11.09&14&\parbox[t]{1.5cm}{\{Blijdorpplein, Kleinpolderplein\}}&\parbox[t]{3cm}{\{Kerkhoflaan Buffer, Hofwijk, return-West-Sidelinge, Delft, Vrijenban, 2e Hogenbanweg\}}\\
        0.4 &$1.18\times10^6$& 500000
&+41\%&16.61&8.98&13.59&16&\parbox[t]{1.5cm}{\{Blijdorpplein\}}&\parbox[t]{3cm}{\{Zaagmolenbrug, return-Nieuwe Crooswijkseweg, Hofwijk, return-West-Sidelinge, Delft, Vrijenban, 2e Hogenbanweg\}}\\
        0.6 &$1.29\times10^6$& 500000
&+55\%&19.11&11.31&15.29&22&\parbox[t]{1.5cm}{\{Blijdorpplein\}}&\parbox[t]{3cm}{\{Zaagmolenbrug, return-Nieuwe Crooswijkseweg, Hofwijk, return-West-Sidelinge, Delft, Vrijenban, 2e Hogenbanweg\}}\\
        0.8 &$1.40\times10^6$& 660000
&+67\%&12.5&13.42&16.82&25&\parbox[t]{1.5cm}{\{return-Gatwickbaan\}}&\parbox[t]{3cm}{\{return-Vliegveldweg, return-Van Der Sasstraat, Blijdorpplein, De Lugt, Delft, Heertjeslaan, return-Nieuwe Crooswijkseweg, Schipluiden-De Zweth, Crooswijksestraat\}}\\
        1 &$1.50\times10^6$& 660000&+79\%&14.72&15.09&18.49&38&\parbox[t]{1.5cm}{\{return-Vliegveldweg\}}&\parbox[t]{3cm}{\{return-Van Der Sasstraat, Blijdorpplein, De Lugt, Delft, Heertjeslaan, return-Nieuwe Crooswijkseweg, Schipluiden-De Zweth, Crooswijksestraat, return-Gatwickbaan\}}\\
    \bottomrule
    \end{tabular}}\\
{\textbf{OFV}: Objective function value or total costs; \textbf{CS installation}: Charging station installation; \textbf{SS}: standard charging station; \textbf{FF}: fast-feeding charging station; Network configuration contains stop names (shown in \Cref{fig:bus-network}) that specific charging station types is installed}
    \label{tab:BoU-performance}
\end{table}

\begin{table}[htbp]
    \centering
    \renewcommand{\arraystretch}{1.5} 
    \captionsetup{justification=raggedright,singlelinecheck=false,skip=2pt} 
    \caption{\textcolor{blue}{DRCC-CID solution sensitivity in terms of budget parameter ($\theta$)}}
    \tiny
    \textcolor{blue}{
    \begin{tabular}{lp{1.5cm}p{1.125cm}p{1.5cm}p{1cm}p{1cm}p{1cm}p{0.5cm}p{1.5cm}p{2.5cm}}
    \toprule
       \multirow{2}{*}{\textbf{$\theta$}}& \multicolumn{3}{l}{\textbf{Electrification cost}}&\multicolumn{3}{l}{\textbf{Battery capacity (kWh)}} & \multirow{2}{*}{\textbf{CPU}(s)}&\multicolumn{2}{l}{\textbf{Network configuration}}\\
       \cmidrule{2-7}
       \cmidrule{9-10}
        &OFV (\euro)& CS installation&Relative difference&Line 33 &Line 38 &Line 40& &  SS&FF\\ 
       \hline
        0 & $8.413\times10^5$&  340000&\textbf{benchmark}.&9.26&7.5&11.87&9.41&\parbox[t]{1.5cm}{\{Delft, Balthasar Van Der Polweg\} }&\parbox[t]{3cm}{\{Paradijsplein, Schipluiden- Halfwege, return-De Lugt, 2e Hogenbanweg\}}\\
        0.2
&$9.81\times10^5$& 260000
&17\%&16.61&9.01&15.54&31&\parbox[t]{1.5cm}{\{return-Gatwickbaan\}}&\parbox[t]{3cm}{\{Delft, Ackerdijkseweg, De Lugt, 2e Hogenbanweg, Paradijsplein\}}\\
        0.4
&$9.95\times10^5$& 260000
&18\%&16.55&9.75&15.7&96&\parbox[t]{1.5cm}{\{return-Gatwickbaan\}}&\parbox[t]{3cm}{\{Delft, Ackerdijkseweg, De Lugt, 2e Hogenbanweg, Paradijsplein\}}\\
        0.6
&$9.98\times10^5$& 260000
&19\%&17.14&9.49&15.54&135&\parbox[t]{1.5cm}{\{return-Gatwickbaan\}}&\parbox[t]{3cm}{\{Delft, Ackerdijkseweg, De Lugt, 2e Hogenbanweg, Paradijsplein\}}\\
        0.8
&$9.98\times10^5$& 260000
&19\%&17.12&9.51&15.54&330&\parbox[t]{1.5cm}{\{return-Gatwickbaan\}}&\parbox[t]{3cm}{\{Delft, Ackerdijkseweg, De Lugt, 2e Hogenbanweg, Paradijsplein\}}\\
        1&$1.005\times10^6$& 260000
&20\%&17.02&10.05&15.54&393&\parbox[t]{1.5cm}{\{return-Gatwickbaan\}}&\parbox[t]{3cm}{\{Delft, Ackerdijkseweg, De Lugt, 2e Hogenbanweg, Paradijsplein\}}\\
        2.5&$1.005\times10^6$& 260000&20\%&17.02&10.05&15.54&1296&\parbox[t]{1.5cm}{\{return-Gatwickbaan\}}&\parbox[t]{3cm}{\{Delft, Ackerdijkseweg, De Lugt, 2e Hogenbanweg, Paradijsplein\}}\\
    \bottomrule
    \end{tabular}\\
{\textbf{OFV}: Objective function value or total costs; \textbf{CS installation}: Charging station installation; \textbf{SS}: standard charging station; \textbf{FF}: fast-feeding charging station; Network configuration contains stop names  (shown in \Cref{fig:bus-network}) that specific charging station types is installed}
    \label{tab:DRCC-performance}}
\end{table}

\newpage
\section{Simulation pseudo-code}
\label{simulation_pseudo}
\begin{algorithm}
\caption{Simulation Pseudo-Code for Feasibility of CID Problem}
\textcolor{blue}{\begin{algorithmic}[1]
\STATE \textbf{Input:} Charging station locations and types ($x_{st}$), battery capacity ($z_k$), number of simulation scenarios ($N$), energy consumption distribution
\STATE \textbf{Output:} Feasibility results and violations per scenario for each bus line
\FOR{each bus line $k$}
    \STATE Initialize feasibility flag: $\text{feasible} \leftarrow \text{True}$
    \FOR{scenario $s = 1$ to $N$}
        \STATE Calculate initial battery capacity at depot: $B_k \leftarrow \bar{b} z_k$
        \STATE Set initial battery level: $b_k \leftarrow B_k$
        \STATE Generate random energy consumption $EC_{(i,j)}$ between consecutive stops for line $k$
        \FOR{each consecutive stops $(i, j)$ on line $k$}
            \IF{Charging station is installed at $i$: $x_{it} = 1$}
                \STATE Calculate charging received: $\text{Charging received} \leftarrow P_t \Delta_{ks}$
                \STATE Update battery level: $b_k \leftarrow \min(B_k, b_k + \text{Charging received})$
            \ENDIF
            \STATE Update battery level after consumption: $b_k \leftarrow b_k - EC_{(i,j)}$
            \IF{$b_k < \underline{b}z_k$}
                \STATE Set $\text{feasible} \leftarrow \text{False}$ and record violation
                \STATE \textbf{break}
            \ENDIF
        \ENDFOR
    \ENDFOR
\ENDFOR
\STATE Compute feasibility rate of bus line $k$: $\frac{\text{Number of Feasible Simulations}}{N} \times 100\%$
\end{algorithmic}}
\end{algorithm}

\newpage

\section{DRCC preformance sample size analysis}
\label{appendix:quantity}

\begin{table}[h]
\centering
    \renewcommand{\arraystretch}{1.5} 
    \captionsetup{justification=raggedright,singlelinecheck=false,skip=2pt} 
    \centering

    \caption{\textcolor{blue}{Impact of Sample Size on the Feasibility of DRCC-CID Solutions}}
\tiny
\label{tab:DRCC-quantity}
\textcolor{blue}{
    \begin{tabular}{lp{1.25cm}p{1.55cm}p{2cm}p{1cm}p{0.5cm}p{0.5cm}p{0.5cm}p{2cm}l}
        \toprule
           nr.&\textbf{Design distribution}&\textbf{Out-of-sampling distribution}&\textbf{Uncertainty level of DRCC-CID}& \textbf{Data sample size}&\multicolumn{3}{l}{\textbf{Bus line number}}&\textbf{Average feasibility of the bus network}&\textbf{CPU} (s)\\
           \cmidrule{6-8}&&& &  &38& 33&40& &\\
        \midrule
           1&$\mathcal{U}[0, \hat{\mu}]$&$\mathcal{U}[0, \textbf{1.2}\times \hat{\mu}]$& low: ($\theta = 0.2$)& 100&14\%& 48\%&49\%&37\% &24\\
           &&& & 300&39\%& 48\%&43\%&43\% &164\\
            &&& & 500&47\%& 48\%&58\%&51\% &853\\
            \cmidrule{4-9}&&& high: ($\theta = 0.8$)& 100&53\%& 50\%&65\%&56\% &428\\
            &&& & 300&56\%& 64\%&75\%&65\% &799\\
 & & & & 500& 57\%& 77\%& 73\%&69\% &2923\\
        \midrule
          2&$\mathcal{T}[0, \hat{\mu}, \frac{\hat{\mu}}{2}]$&$\mathcal{T}[0, \hat{\mu}, \hat{\mu}]$& low: ($\theta = 0.2$)& 100&4\%& 20\%&20\%&15\% &19\\
          &&& & 300&4\%& 6\%&4\%&5\% &193\\
            &&& & 500&17\%& 12\%&17\%&15\% &118\\
         \cmidrule{4-9}&&& high: ($\theta = 0.8$)& 100&17\%& 20\%&13\%&17\% &101\\
            &&& & 300&8\%& 10\%&19\%&12\% &799\\
 & & & & 500& 11\%& 29\%& 46\%&29\% &5231\\
 \bottomrule
    \end{tabular}} 
\end{table}

\end{document}